\documentclass[leqno]{article}
\usepackage{amssymb}
\usepackage{color}
 \def\qed{\unskip\quad \hbox{\vrule\vbox to 6pt {\hrule width
          4pt\vfill\hrule}\vrule} }

                         \newtheorem{Th}{Theorem}
                         \newtheorem{Prop}{Proposition}
                         \newtheorem{Lemma}{Lemma}

                         \newtheorem{Def}{Definition}

\newtheorem{question}{Question}

\newcommand{\fre}{\rule[-0.08cm]{0.06cm}{0.38cm}}

\newcommand{\tor}{{\Bbb{T}}}

\newcommand{\cs}{{\cal S}}

\newcommand{\cc}{{\cal C}}
\newcommand{\ct}{{\cal T}}
\newcommand{\cu}{{\cal U}}

\newcommand{\cb}{{\cal B}}
\newcommand{\ca}{{\cal A}}

\newcommand{\cp}{{\cal P}}

\newcommand{\crr}{{\cal R}}
\newcommand{\lf}{\Lambda_\varphi}
\newcommand{\svf}{\Sigma_\varphi}

\newcommand{\A}{\mathbb{A}}
\newcommand{\Q}{\mathbb{Q}}
\newcommand{\R}{{\mathbb{R}}}
\newcommand{\T}{{\mathbb{T}}}
\newcommand{\C}{{\mathbb{C}}}
\newcommand{\Z}{{\mathbb{Z}}}
\newcommand{\N}{{\mathbb{N}}}
\newcommand{\xbm}{(X,{\cal B},\mu)}

\newcommand{\ycn}{(Y,{\cal C},\nu)}

\newcommand{\stfs}{\sigma_{T_{\varphi,{\cal S}}}}
\newcommand{\tfs}{T_{\varphi,{\cal S}}}

\newcommand{\ov}{\overline}
\newcommand{\beq}{\begin{equation}}
\newcommand{\eeq}{\end{equation}}
\newcommand{\vep}{\varepsilon}
\newcommand{\va}{\varphi}
\newcommand{\si}{\sigma}
\newcommand{\ot}{\otimes}
\newcommand{\la}{\lambda}

\begin{document}
\title{Spectral theory of dynamical systems}
\author{Adam Kanigowski\footnote{Department of Mathematics,  University of Maryland at College Park, College Park,MD
20740, USA, adkanigowski@gmail.com} \ and \  Mariusz Lema\'nczyk\footnote{Faculty of Mathematics and Computer Science,
Nicolaus Copernicus University, 87-100 Toru\'n, Chopin street
12/18, Poland, mlem@mat.umk.pl}}
\date{}

\maketitle

\scriptsize

\tableofcontents

\thispagestyle{empty}



\normalsize
\section*{Glossary and notation}
\addcontentsline{toc}{section}{Glossary and notation} {\bf
Spectral decomposition of a unitary representation.} If
$\cu=(U_a)_{a\in\A}$ is a continuous unitary representation of a
locally compact second countable (l.c.s.c.) Abelian group $\A$ in
a separable Hilbert space $H$ then a decomposition
$H=\bigoplus_{i=1}^\infty \A(x_i)$ is called {\em spectral} if
$\sigma_{x_1}\gg\sigma_{x_2}\gg\ldots$ (such  a sequence of
measures is also called {\em spectral}); here
$\A(x):=\ov{\rm span}\{U_ax:\:a\in\A\}$ is called the {\em cyclic
space} generated by $x\in H$ and $\sigma_x$ stands for the
spectral measure of $x$.

\vspace{1ex}

{\bf Maximal spectral type and the multiplicity function of
$\cu$.} \ The {\em maximal spectral type} $\sigma_{\cu}$ of $\cu$
is the type of $\sigma_{x_1}$ in any spectral decomposition of
$H$; the {\em multiplicity function}
$M_{\cu}:\widehat{\A}\to\{1,2,\ldots\}\cup\{+\infty\}$ is defined
$\sigma_{\cu}$-a.e. and $ M_{\cu}(\chi)=\sum_{i=1}^\infty
1_{Y_i}(\chi)$, where $Y_1=\widehat{\A}$ and
$Y_i={\rm supp}\,\frac{d\sigma_{x_i}}{d\sigma_{x_1}}$ for $i\geq2$.

A representation $\cu$ is said to have {\em simple spectrum} if
$H$ is reduced to a single cyclic space. The multiplicity is {\em
uniform} if there is only one essential value of $M_{\cu}$. The
essential supremum of $M_{\cu}$ is called the {\em maximal
spectral multiplicity}.  $\cu$ is said to have {\em discrete
spectrum} if $H$  has an orthonormal basis consisting of
eigenvectors of $\cu$; $\cu$ has {\em singular} ({\em Haar,
absolutely continuous}) spectrum if the maximal spectral type of
$\cu$ is singular with respect to (equivalent to, absolutely continuous with) a Haar measure of $\widehat{\A}$.

\vspace{1ex}

{\bf Koopman representation of a dynamical system $\ct$.} \ Let
$\A$ be a l.c.s.c. (and not compact) Abelian group and
$\ct:a\mapsto T_a$ a representation of $\A$ in the group $Aut\xbm$
of (measure-preserving) automorphisms of a standard probability
Borel space $\xbm$. The {\em Koopman representation}
$\cu=\cu_{\ct}$ of $\ct$ in $L^2\xbm$ is defined as the unitary
representation $a\mapsto U_{T_a}\in U( L^2\xbm)$, where
$U_{T_a}(f)=f\circ T_a$.

\vspace{1ex}

{\bf Ergodicity, weak mixing, mild mixing, mixing and rigidity of
$\ct$.} \ A measure-preserving $\A$-action $\ct=(T_a)_{a\in\A}$ is
called {\em ergodic} if $\chi_0\equiv1\in\widehat{\A}$ is a simple
eigenvalue of $\cu_{\ct}$. It is {\em weakly mixing} if
$\cu_{\ct}$ has a continuous spectrum on the subspace $L^2_0\xbm$
of zero mean functions. $\ct$ is said to be {\em rigid} if there
is a sequence $(a_n)$ going to infinity in $\A$ such that the
sequence  $(U_{T_{a_n}})$ goes  to the identity in the strong (or
weak) operator topology; $\ct$ is said to be {\em mildly mixing}
if it has no non-trivial rigid factors. We say that $\ct$ is {\em
mixing} if the operator equal to zero is the only limit point of
$\{U_{T_a}|_{L^2_0\xbm}:\:a\in\A\}$ in the weak operator topology.

\vspace{1ex}

{\bf Spectral disjointness.} \ Two $\A$-actions $\cs$ and $\ct$
are called {\em spectrally disjoint} if the maximal spectral types
of their Koopman representations $\cu_{\ct}$ and $\cu_{\cs}$ on
the corresponding $L^2_0$-spaces are mutually singular.

\vspace{1ex}

{\bf SCS property}. We say that a Borel measure $\sigma$ on
$\widehat{\A}$ satisfies the {\em strong convolution singularity}
property (SCS property) if, for each $n\geq1$, in the
disintegration (given by the map
$(\chi_1,\ldots,\chi_n)\mapsto\chi_1\cdot\ldots\cdot\chi_n$)
$\sigma^{\otimes
n}=\int_{\widehat{\A}}\nu_\chi\,d\sigma^{(n)}(\chi)$ the
conditional measures $\nu_\chi$ are atomic with exactly $n!$ atoms
($\sigma^{(n)}$ stands for the $n$-th convolution
$\sigma\ast\ldots\ast\sigma$). An $\A$-action $\ct$ satisfies the
SCS property if the maximal spectral type of $\cu_{\ct}$ on
$L^2_0$ is a type of an SCS measure.

\vspace{1ex}

{\bf Kolmogorov group property.} \ An $\A$-action $\ct$ satisfies
the {\em Kolmogorov group property} if
$\sigma_{\cu_{\ct}}\ast\sigma_{\cu_{\ct}}\ll\sigma_{\cu_{\ct}}$.

\vspace{1ex}

{\bf Weighted operator.} \ Let $T$ be an ergodic automorphism of
$\xbm$  and $\xi:X\to\T$ be a measurable function. The (unitary)
operator $V=V_{\xi,T}$ acting on $L^2\xbm$ by the formula
$V_{\xi,T}(f)(x)=\xi(x)f(Tx)$ is called a {\em weighted operator}.

\vspace{1ex}

{\bf Induced automorphism.} \ Assume that $T$ is an
automorphism of a standard probability Borel space $\xbm$. Let
$A\in\cb$, $\mu(A)>0$. The {\em induced automorphism} $T_A$ is
defined on the conditional space $(A,\cb_A,\mu_A)$, where $\cb_A$
is the trace of $\cb$ on $A$, $\mu_A(B)=\mu(B)/\mu(A)$ for
$B\in\cb_A$ and $T_A(x)=T^{k_A(x)}x$, where $k_A(x)$ is the
smallest $k\geq1$ for which $T^kx\in A$.

\vspace{1ex}

{\bf AT property of an automorphism.} \ An automorphism $T$ of a
standard  probability Borel space $\xbm$ is called {\em
approximatively transitive} (AT for short) if for every $\vep>0$
and every finite set $f_1,\ldots,f_n$ of non-negative
$L^1$-functions on $\xbm$ we can find $f\in L^1\xbm$ also
non-negative such that $ \|f_i-\sum_{j}\alpha_{ij}f\circ
T^{n_j}\|_{L_1}<\vep$ for all $i=1,\ldots, n$ (for some
$\alpha_{ij}\geq 0$, $n_j\in\N$).

\vspace{1ex}

{\bf Cocycles and group extensions.} \ If $T$ is an ergodic
automorphism, $G$ is a l.c.s.c.\ Abelian group and $\va:X\to G$ is
measurable then the pair $(T,\va)$ generates a {\em cocycle}
$\va^{(\cdot)}(\cdot):\Z\times X\to G$, where
$$
\va^{(n)}(x)=\left\{\begin{array}{lll}
\va(x)+\ldots+\va(T^{n-1}x)&
\mbox{for}&n>0,\\
0&\mbox{for}&n=0\\
-(\va(T^nx)+\ldots+\va(T^{-1}x))&\mbox{for}&n<0.\end{array}\right.
$$
(That is, $(\va^{(n)})$ is a standard  1-cocycle in the algebraic
sense for the $\Z$-action $n(f)=f\circ T^n$ on the group of
measurable functions on $X$ with values in $G$; hence the function
$\va:X\to G$ itself is often called a {\em cocycle}.)

Assume additionally that $G$ is compact. Using the cocycle $\va$,
we define a {\em group extension} $T_{\va}$  on $(X\times
G,\cb\ot\cb(G),\mu\ot \la_G)$ ($\la_G$ stands for Haar measure of
$G$), where $T_{\va}(x,g)=(Tx,\va(x)+g)$.

\vspace{1ex}

{\bf Special flow.} \ Given an ergodic automorphism $T$ on a
standard probability Borel space $\xbm$ and a positive integrable
function $f:X\to\R^+$ we put
$$ X^f=\{(x,t)\in X\times\R:\:0\leq
t<f(x)\},\;\cb^f=\cb\ot\cb(\R)|_{X^f},$$ and define $\mu^f$  as
normalized $\mu\otimes\lambda_{\R}|_{X^f}$. By a {\em special
flow} we mean the $\R$-action $T^f=(T^f_t)_{t\in\R}$ under which a
point $(x,s)\in X^f$ moves vertically with the unit speed, and
once it reaches the graph of $f$, it is identified with $(Tx,0)$.
\vspace{1ex}

{\bf Time change.} \  Let $\mathcal{R}=(R_t)_{t\in\R}$ be a flow on $\xbm$  and let $v\in L^1\xbm$ be a positive function. The function $v$ determines a cocycle over $\mathcal{R}$  given by the formula
$$
v(t,x):= \int_0^t v(R_sx) d s.
$$
Then for a.e.\ $x\in {X}$ and all $t\in\R$, there exists a unique $u=u(t,x)$  such that
$$
\int_0^uv(R_sx)\, ds=t.$$
 Now, we can define the flow $\widetilde{R}_t(x):=R_{u(t,x)}(x)$. The new flow $\widetilde{\mathcal{R}}=(\widetilde{R}_t)_{t\in\R}$ has the same orbits as the original flow, and it preserves the measure $\widetilde{\mu}\ll\mu$ (hence it is ergodic if $\mathcal{R}$ was), where $\frac{d\widetilde{\mu}}{d\mu}=v/\int_Xv\,d\mu$.

\vspace{1ex}

{\bf Markov operator.} \
 A linear operator $J:L^2\xbm\to L^2\ycn$ is called
{\em Markov} if it sends non-negative functions to non-negative
functions and $J1=J^\ast1=1$.

\vspace{1ex}

{\bf Unitary actions on Fock spaces.} \ If $H$ is a separable
Hilbert space then by $H^{\odot n}$ we denote the subspace of
$n$-tensors of $H^{\otimes n}$ symmetric under all permutations
of coordinates, $n\geq1$; then the Hilbert space
$F(H):=\bigoplus_{n=0}^\infty H^{\odot n}$ is called a {\em
symmetric Fock space}.  If $V\in U(H)$ then $F(V):=\bigoplus_
{n=0}^\infty V^{\odot n}\in U(F(H))$, where $ V^{\odot n}=
V^{\otimes n}|_{ H^{\odot n}}$.

\section{Definition of the subject}
Spectral theory of dynamical systems is a study of special unitary
representations, called Koopman representations (see the
glossary). Invariants of such representations are called spectral
invariants of measure-preserving systems. Together with the
entropy, they consitute the most important invariants used in the
study of measure-theoretic intrinsic properties and classification
problems of dynamical systems as well as in applications. Spectral
theory was originated by von Neumann, Halmos and Koopman in the
1930s. In this article we will focus on recent progresses in the
spectral theory of finite measure-preserving dynamical systems.

\section{Introduction}\label{czesc0}
 Throughout $\A$ denotes a non-compact l.c.s.c.\
Abelian group ($\A$ will be most often $\Z$ or $\R$). The
assumption of second countability implies that $\A$ is metrizable,
$\sigma$-compact and the space $C_0(\A)$ is separable. Moreover
the dual group $\widehat{\A}$ is also l.c.s.c. Abelian.
\subsection{General unitary representations}
 We are interested in {\em unitary}, that is with values in the unitary
group $U(H)$ of a Hilbert space $H$, (weakly) continuous
representations $V: \A\ni a\mapsto V_a\in U(H)$ of such groups
(the scalar valued maps $a\mapsto\langle V_ax,y\rangle$ are
continuous for each $x,y\in H$).

 Let
$H=L^2(\widehat{\A},\cb(\widehat{\A}),\mu)$, where
$\cb(\widehat{\A})$ stands for the $\sigma$-algebra of Borel sets
of $\widehat{\A}$ and $\mu\in M^+(\widehat{\A})$ (whenever $X$ is
a l.c.s.c. space, by $M(X)$ we denote the set of complex Borel
measures on $X$, while $M^+(X)$ stands for the subset of positive
(finite) measures). Given $a\in\A$, for $f\in
L^2(\widehat{\A},\cb(\widehat{\A}),\mu)$ put
$$
V^{\mu}_a(f)(\chi)=i(a)(\chi)\cdot f(\chi)=\chi(a)\cdot
f(\chi)\;\;(\chi\in\widehat{\A}),$$ where
$i:\A\to\widehat{\widehat{\A}}$ is the canonical Pontriagin
isomorphism of $\A$ with its second dual. Then
$V^{\mu}=(V^{\mu}_a)_{a\in\A}$ is a unitary representation of
$\A$. Since $C_0(\widehat{\A})$ is dense in
$L^2(\widehat{\A},\mu)$, the latter space is separable. Therefore
also direct sums $\bigoplus_{i=1}^\infty V^{\mu_i}$ of such type
representations will be unitary representations of $\A$ in
separable Hilbert spaces.

\begin{Lemma}[Wiener Lemma] If $F\subset L^2(\widehat{\A},\mu)$ is
a closed $V^\mu_a$-invariant subspace for all $a\in\A$ then
$F=1_{Y}L^2(\widehat{\A},\cb(\widehat{\A}),\mu)$ for some Borel
subset $Y\subset \widehat{\A}$.
\end{Lemma}

Notice however that since $\A$ is not compact (equivalently,
$\widehat{\A}$ is not discrete), we can find $\mu$ continuous and
therefore $V^\mu$ has no irreducible (non-zero) subrepresentation.
From now on only unitary representations of $\A$ in separable
Hilbert spaces will be considered and we will show how to classify
them.

A function $r:\A\to\C$ is called {\em positive definite} if
\beq\label{dodokr} \sum_{n,m=0}^Nr(a_n-a_m)z_n\ov{z_m}\geq0\eeq
for each $N>0$, $(a_n)\subset \A$ and $(z_n)\subset \C$. The
central result about positive definite functions is the following
theorem (see e.g. \cite{Rudin}).

\begin{Th}[Bochner-Herglotz] Let $r:\A\to\C$ be continuous. Then
$r$ is positive definite if and only if there exists (a unique)
$\sigma\in M^+(\widehat{\A})$ such that
$$
r(a)=\int_{\widehat{\A}}\chi(a)\,d\sigma(\chi)\;\;\mbox{for
each}\;a\in\A.$$
\end{Th}

If now $\cu=(U_a)_{a\in\A}$ is a representation of $\A$ in $H$
then for each $x\in H$ the function $r(a):=\langle U_ax,x\rangle$
is continuous and satisfies~(\ref{dodokr}),  so it is positive
definite. By the Bochner-Herglotz Theorem there exists a unique
measure $\sigma_{\cu,x}=\sigma_x\in M^+(\widehat{\A})$ (called the
{\em spectral measure} of $x$) such that
$$
\widehat{\sigma}_x(a):=\int_{\widehat{\A}}i(a)(\chi)\,d\sigma_x(\chi)=\langle
U_ax,x\rangle$$ for each $a\in\A$.  Since the partial map
$U_ax\mapsto i(a)\in L^2(\widehat{\A},\sigma_x)$ is isometric and
equivariant, there exists a unique extension of it to a unitary
operator $W:\A(x)\to L^2(\widehat{\A},\sigma_x)$ giving rise to an
isomorphism of $\cu|_{\A(x)}$ and $V^{\sigma_x}$. Then the
existence of a spectral decomposition is proved by making use of
separability and a choice of maximal cyclic spaces at every step
of an induction procedure. Moreover, a spectral decomposition is
unique in the following sense.

\begin{Th}[Spectral Theorem]\label{twsp}
If $H=\bigoplus_{i=1}^\infty \A(x_i)=\bigoplus_{i=1}^\infty
\A(x'_i)$ are two spectral decompositions of $H$ then $
\sigma_{x_i}\equiv\sigma_{x'_i}\;\;\mbox{for each}\;\;i\geq1. $
\end{Th}

It follows that the representation $\cu$ is entirely determined by
the types (the sets of equivalent measures to a given one) of a
decreasing sequence of measures or,  equivalently, $\cu$ is
determined by its maximal spectral type $\sigma_{\cu}$ and its
multiplicity function $M_{\cu}$.

 Notice that eigenvalues of $\cu$
correspond to Dirac measures: {\em $\chi\in\widehat{\A}$ is an
eigenvalue (i.e.\ for some $\|x\|=1$, $U_a(x)=\chi(a)x$ for each
$a\in\A$) if and only if $\sigma_{\cu,x}=\delta_{\chi}$.}
Therefore $\cu$ has a discrete spectrum if and only if the maximal
spectral type of $\cu$ is a discrete measure.

 We refer the reader to \cite{Gl}, \cite{Ka-Th},
\cite{Le}, \cite{Na}, \cite{Pa} for presentations of spectral
theory needed in the theory of dynamical systems -- such
presentations are usually given for $\A=\Z$ but once we have the
Bochner-Herglotz Theorem and the Wiener Lemma, their extensions to
the general case are straightforward.

\subsection{Koopman representations}
We will consider {\em measure-preserving} representations of $\A$.
It means that we fix a probability standard Borel space $\xbm$ and
by $Aut\xbm$ we denote the group of automorphisms of this space,
that is, $T\in Aut\xbm$ if $T:X\to X$ is a bimeasurable (a.e.)
bijection satisfying $\mu(A)=\mu(TA)=\mu(T^{-1}A)$ for each
$A\in\cb$. Consider then a representation of $\A$ in $Aut\xbm$,
that is, a group homomorphism $a\mapsto T_a\in Aut\xbm$; we write
$\ct=(T_a)_{a\in\A}$. Moreover, we require that the associated
Koopman representation $\cu_{\ct}$ is continuous. Unless
explicitly stated, $\A$-actions are assumed to be {\em free}, that
is, we assume that for $\mu$-a.e.\ $x\in X$ the map $a\mapsto T_ax$
is 1-1.  In fact, since constant functions are obviously invariant
for $U_{T_a}$, equivalently, the trivial character~1 is always an
eigenvalue of $\cu_{\ct}$, the Koopman representation is
considered only on the subspace $L^2_0\xbm$ of zero mean
functions.  We will restrict our attention only to {\em ergodic}
dynamical systems (see the glossary). It is easy to see that $\ct$
is ergodic if and only if  whenever $A\in\cb$ and $A=T_a(A)$
($\mu$-a.e.) for all $a\in A$ then $\mu(A)$ equals~$0$ or~$1$.  In
case of ergodic Koopman representations, all eigenvalues are
simple. In particular, (ergodic) Koopman representations with
discrete spectra have simple spectra. The reader is referred to
monographs mentioned above as well as to \cite{Co-Fo-Si},
\cite{Pe}, \cite{Ru}, \cite{Si}, \cite{Wa} for basic facts on the
theory of dynamical systems. See also survey articles \cite{LeSurvey} and \cite{DaSurvey}.

The passage $\ct\mapsto\cu_{\ct}$ can be seen as functorial
(contravariant). In particular, a measure-theoretic isomorphism of
$\A$-systems $\ct$ and $\ct'$ implies spectral isomorphism of the
corresponding Koopman representations; hence  spectral properties
are measure-theoretic invariants. Since unitary representations
are completely classified, one of the main questions in the
spectral theory of dynamical systems is to decide which pairs
$([\sigma],M)$ can be realized by Koopman representations. The
most spectacular is the Banach problem concerning a realization, for $\A=\Z$, of
$([\la_{\T}],M\equiv1)$, see Section~\ref{s:Banach}. Another important problem  is to give a
complete spectral classification in some  classes of dynamical
systems (classically, it was done in the theory of Kolmogorov and
Gaussian dynamical systems). We will also see how spectral
properties of dynamical systems can determine their statistical
(ergodic) properties; a culmination given  by the fact that a
spectral isomorphism may imply measure-theoretic similitude
(discrete spectrum case, Gaussian-Kronecker case). An old conjecture is
that a dynamical system $\ct$ either is spectrally determined or
there are uncountably many pairwise non-isomorphic systems
spectrally isomorphic to $\ct$.

We could also consider Koopman representations in $L^p$ for $1\leq
p\neq2$. However, whenever $W:L^p\xbm\to L^p\ycn$ is a surjective
isometry and $W\circ U_{T_a}=U_{S_a}\circ W$ for each $a\in\A$
then by the Lamperti Theorem (e.g.\ \cite{Ro}) the isometry $W$
has to come from a non-singular pointwise map $R:Y\to X$ and, by
ergodicity,  $R$ ``preserves" the measure $\nu$ and hence
establishes a measure-theoretic isomorphism \cite{rusek} (see also
\cite{Le}). Thus spectral classification of such Koopman
representations goes back to the measure-theoretic classification
of dynamical systems, so it looks hardly interesting. This does
not mean that there are no interesting dynamical questions for
$p\neq2$. Let us mention still open Thouvenot's question
(formulated in the 1980s) for $\Z$-actions: {\em  For every
ergodic $T$ acting on $\xbm$, does there exist $f\in L^1\xbm$ such
that the closed linear span of $f\circ T^n$, $n\in \Z$, equals
$L^1\xbm$?}

Iwanik  \cite{Iw1}, \cite{Iw2} proved that if $T$ is a  system
with positive entropy then its  $L^p$-multiplicity is~$\infty$ for
all $p>1$. Moreover, Iwanik and  de Sam Lazaro  \cite{Iw-SamL}
proved that for Gaussian systems (they will be considered in
Section~\ref{czesc10}) the $L^p$--multiplicities are the same for
all $p>1$ (see also \cite{Le-Sam de Lazaro}).

\subsection{Markov operators, joinings and Koopman representations,
disjointness and spectral disjointness, entropy} We would like to
emphasize that spectral theory is closely related to the theory of
joinings  (see de la Rue's article \cite{Thierry} for needed
definitions). The elements $\rho$ of the set $J(\cs,\ct)$ of
joinings of two $\A$-actions $\cs$ and $\ct$ are in a 1-1
correspondence with Markov operators $J=J_\rho$ between the
$L^2$-spaces equivariant with the corresponding Koopman
representations (see the glossary and \cite{Thierry}). The set of
ergodic self-joinings of an ergodic $\A$-action  $\ct$ is denoted
by  $J^e_2(\ct)$.

Each Koopman representation $\cu_{\ct}$ consists of Markov
operators (indeed, $U_{T_a}$ is clearly a Markov operator). In
fact, if $U\in U(L^2\xbm)$ is Markov then it is of the form $U_R$,
where $R\in Aut\xbm$ (\cite{Le-Pa2}). This allows us to see
Koopman representations as unitary Markov representations, but
since a spectral isomorphism does not ``preserve" the set of
Markov operators, spectrally isomorphic systems can have
drastically different sets of self-joinings.

We will touch here only some aspects of interactions (clearly, far
from completeness) between the spectral theory and the theory of
joinings.


In order to see however an example of such interactions  let us
recall that the simplicity of eigenvalues for ergodic systems
yields a short ``joining" proof  of the classical isomorphism
theorem of Halmos-von Neumann in the discrete spectrum case: {\em
Assume that $\cs=(S_a)_{a\in\A}$ and $\ct=(T_a)_{a\in\A}$ are
ergodic $\A$-actions on $\xbm$ and $\ycn$ respectively. Assume
that both Koopman representations have purely discrete spectrum
and that their group of eigenvalues are the same. Then $\cs$ and
$\ct$ are measure-theoretically isomorphic.} Indeed, each ergodic
joining of $\ct$ and $\cs$ is the graph of an isomorphism of these
two systems (see \cite{Le}; see also Goodson's
proof in \cite{Go}). Another example of such interactions appear
in the study Rokhlin's multiple mixing problem and its relation
with the {\em pairwise independence property} (PID) for joinings
of higher order. We will not deal with this subject here,
referring the reader to \cite{Thierry} (see however
Section~\ref{liftingmix}).

Following \cite{Fu1}, two  $\A$-actions $\cs$ and $\ct$ are
called {\em disjoint} provided the product measure is the only
element in $J(\cs,\ct)$ (if they are disjoint, one of these actions has to be ergodic).  It was already noticed in \cite{Ha-Pa}
that spectrally disjoint systems are disjoint in the Furstenberg
sense; indeed, ${\rm Im}(J_\rho|_{L^2_0})=\{0\}$ since
$\sigma_{\ct,J_\rho f}\ll\sigma_{\cs,f}$.

Notice that whenever we take $\rho\in J^e_2(\ct)$ we obtain a new
ergodic $\A$-action $(T_a\times T_a)_{a\in\A}$ defined on the
probability space $(X\times X,\rho)$. One can now ask how close
spectrally to $\ct$  is this new action? It turns out that except
of the obvious fact that the marginal $\sigma$-algebras are
factors, $(\ct\times\ct,\rho)$ can have other factors spectrally
disjoint from $\ct$:  the most striking phenomenon here is a
result of Smorodinsky and Thouvenot \cite{Sm-Th} (see also
\cite{Da-Park}) saying that each zero entropy system is a factor
of an ergodic self-joining system of a fixed Bernoulli system
(Bernoulli systems themselves have countable Haar spectrum). The
situation changes if $\rho=\mu\otimes\mu$. In this case for
$f,g\in L^2\xbm$ the spectral measure of $f\otimes g$ is equal to
$\sigma_{\ct,f}\ast\sigma_{\ct,g}$. A consequence of this
observation is that an ergodic action $\ct=(T_a)_{a\in\A}$ is
weakly mixing (see the glossary)  if and only if  the product
measure $\mu\ot\mu$ is an ergodic self-joining of $\ct$.

The entropy which is a basic measure-theoretic invariant  does not
appear when we deal with spectral properties. We will not give
here any formal definition of entropy for amenable group actions
referring the reader to \cite{Or-We}. Assume that $\A$ is
countable and discrete. We always assume that $\A$ is Abelian,
hence it is amenable. For each dynamical system
$\ct=(T_a)_{a\in\A}$ acting on $\xbm$, we can find a largest
invariant sub-$\sigma$ field $\cp\subset\cb$, called the {\em
Pinsker} $\sigma$-algebra, such that the entropy of the
corresponding quotient system  is zero. Generalizing the classical
Rokhlin-Sinai Theorem (see also \cite{Kam} for $\Z^d$-actions),
Thouvenot (unpublished) and independently Dooley and Golodets
\cite{Do-Go} proved this theorem for groups even more general than
those considered here: {\em The spectrum of $\cu_{\ct}$ on
$L^2\xbm\ominus L^2(\cp)$ is Haar with uniform infinite
multiplicity.} This general result is quite intricate and based on
methods introduced to entropy theory by Rudolph and Weiss
\cite{Ru-We} with a very surprising use of Dye's Theorem on
orbital equivalence of all ergodic systems. For $\A$  which is not
countable the same result was proved  in
\cite{uczenWeissa} in case of unimodular amenable groups which are
not increasing union of compact subgroups. It follows that
spectral theory of dynamical systems  essentially reduces to the
zero entropy case.

\section{Maximal spectral type of a Koopman
representation, Alexeyev's Theorem}\label{czesc4} Only few general
properties of maximal spectral types of Koopman representations
are known. The fact that a Koopman representation preserves the
space of real functions implies that its maximal spectral type is
the type of a symmetric (invariant under the map
$\chi\mapsto\ov{\chi}$) measure.

Recall that the {\em Gelfand spectrum} $\sigma(\cu)$ of a
representation $\cu=(U_a)_{a\in\A}$ is defined as the set of {\em
approximative eigenvalues} of $\cu$, that is,
$\sigma(\cu)\ni\chi\in\widehat{\A}$ if for a sequence $(x_n)$
bounded and bounded away from zero, $\|U_ax_n-\chi(a)x_n\|\to 0$
for each $a\in\A$.  The spectrum is a closed subset in the
topology of pointwise convergence, hence, in the
compact-open topology of $\widehat{\A}$. In case of $\A=\Z$,
the above set $\sigma(U)$ is equal to
$\{z\in\C:\:U-z\cdot Id\;\mbox{is not invertible}\}$.

Assume now that $\A$ is countable and discrete (and Abelian). Then
there exists a F\"olner sequence  $(B_n)_{n\geq1}$ whose elements
tile $\A$ \cite{Or-We}. Take a free and ergodic action
$\ct=(T_a)_{a\in\A}$ on $\xbm$. By \cite{Or-We} for each $\vep>0$
we can find a set $Y_n\in\cb$ such that the sets $T_bY_n$ are
pairwise disjoint for $b\in B_n$ and $\mu(\bigcup_{b\in
B_n}T_bY_n)>1-\vep$. For each $\chi\in\widehat{\A}$, by
considering functions of the form $f_n=\sum_{b\in
B_n}\chi(b)1_{T_bY_n}$, we obtain that $\chi\in\sigma(\cu_{\ct})$.
It follows that the topological support of the maximal spectral
type of the Koopman representation of a free and ergodic action is
full (\cite{Ka-Th}, \cite{Le}, \cite{Na}). The theory of Gaussian
systems shows in particular that there are symmetric measures on
the circle whose topological support is the whole circle but which
cannot be maximal spectral types of Koopman representations.

An open well-known question remains of whether an absolutely continuous
measure $\rho$ is the maximal spectral type of a Koopman
representation if and only if $\rho$ is equivalent to a Haar
measure of $\widehat{\A}$ (this is unknown for $\A=\Z$).

Another interesting question was raised by A.\ Katok
(see \cite{Ka-Le}): {\em Is it true that the topological
supports of all measures in a spectral sequence of a Koopman
representation are full}?  If the answer to this question is
positive then for example the essential supremum of
$M_{\cu_{\ct}}$ is the same on all balls of $\widehat{\A}$.

Notice that the fact that all spectral measures in a spectral
sequence are symmetric means that $\cu_{\ct}$ is isomorphic to
$\cu_{\ct^{-1}}$. A. del Junco \cite{Ju-1} showed that generically
for $\Z$-actions, $T$ and its inverse are not
measure-theoretically isomorphic (in fact, he proved disjointness).

Let $\ct$ be an $\A$-action on $\xbm$. One can ask whether a ``good"
function can realize the maximal spectral type of $\cu_{\ct}$. In
particular, can we find a function $f\in L^\infty\xbm$ that
realizes the maximal spectral type? The answer is given in the
following general theorem (see \cite{Le-Wa}).

\begin{Th}[Alexeyev's Theorem]\label{aleksiejew}
Assume that $\cu=(U_a)_{a\in\A}$ is a unitary representation of
$\A$ in a separable Hilbert space $H$. Assume that $F\subset H$ is
a dense linear subspace. Assume moreover that with some $F$-norm
$\fre\cdot\fre$ -- stronger than the norm $\|\cdot\|$ given by the
scalar product --  $F$ becomes a Fr\'echet space. Then, for each
spectral measure $\sigma$ ($\ll\sigma_{\cal U}$) there exists
$y\in F$ such that $\sigma_y\gg\sigma$. In particular, there
exists $y\in F$ that realizes the maximal spectral type.
\end{Th}

By taking $H=L^2\xbm$ and $F=L^\infty\xbm$ we obtain the positive
answer to the original question. Alexeyev \cite{Al} proved the
above theorem for $\A=\Z$ using  analytic functions. Refining
Alexeyev's original proof, Fr\c{a}czek \cite{Fr1} showed the
existence of a sufficiently regular function realizing the maximal
spectral type depending only on the ``regularity" of the
underlying probability space, e.g.\ when $X$ is a compact metric
space (compact smooth manifold) then one can find a continuous (smooth)
function  realizing the maximal spectral type.

By the theory of systems of probabilistic origin (see
Section~\ref{czesc10}), in case of simplicity of the spectrum, one
can easily point out spectral measures whose types are not
realized by (essentially) bounded functions. However, it is still
an open question whether for each Koopman representation
$\cu_{\ct}$ there exists a sequence $(f_i)_{i\geq1}\subset
L^\infty\xbm$ such that the sequence $(\sigma_{f_i})_{i\geq1}$ is
a spectral sequence for $\cu_{\ct}$. For $\A=\Z$ it is unknown
whether the maximal spectral type of a Koopman representation can
be realized by a characteristic function.

The group $Aut\xbm$ considered with the weak operator topology is closed in $U(L^2\xbm)$,\footnote{If we choose $\{A_i:\:i\geq1\}$ a dense subset in $\mathcal{B}$ (considered modulo null sets), then the weak operator topology is metrizable with the metric $d(T_1,T_2):=\sum_{i\geq1}\frac1{2^i}(\mu(T_1(A_i)\triangle T_2(A_i))+\mu(T_1^{-1}(A_i)\triangle T_2^{-1}(A_i)))$.} hence becomes a Polish group. One can then ask what are ``typical'' (largeness is understood as a residual subset) properties of an automorphism of $\xbm$. It is classical (Halmos) that typically an automorphism is weakly mixing,  rigid and has simple spectrum. Some other typical properties will be discussed later. While Halmos already noticed that in the weak operator topology mixing automorphisms form a meager set,  in \cite{Ti1}, S.\ Tikhonov
considers a special (Polish) topology on the set of mixing automorphisms. In fact, this topology was introduced by Alpern \cite{Alp} in 1985 and Tikhonov disproves a conjecture by Alpern by showing that a generic mixing transformation  has simple singular spectrum and is mixing of arbitrary order; moreover, all its powers are disjoint. In \cite{Ti2},  the topology is extended to mixing actions of infinite countable groups $H$, it is given by the metric $d_m$, where for two $H$-actions $\ct_i$ and $H\ni h\mapsto |h|\in \N$ so that $\sum_{h\in H}1/2^{|h|}<+\infty$, we have
$$
d_m(\ct_1,\ct_2):=$$$$\sum_{h\in H}\frac1{2^{|h|}}d(T_{1,h},T_{2,h})+\sup_{h\in H}\sum_{i,j\geq1}\frac1{2^{i+j}}|\mu((T_{1,h}A_i\cap A_j)-\mu(T_{2,h}A_i\cap A_j))|.$$
Bashtanov  \cite{Bash1}, \cite{Bash2} proved that the conjugacy classes (of mixing automorphisms) are dense in this topology. Hence, properties like to have trivial centralizer and no (non-trivial) factors are typical in this topology.

\section{Spectral theory of weighted operators}
\label{czesc5} We now pass to the problem of possible essential
values for the multiplicity function of a Koopman representation.
However, one of known techniques is a use of cocycles, so before
we tackle the multiplicity problem, we will go through some
results concerning spectral theory of compact group extensions
automorphisms which in turn entail a study of weighted operators
(see the glossary).

Assume that $T$ is an ergodic automorphism of a standard Borel
probability space $\xbm$. Let $\xi:X\to\T$ be a measurable
function and let $V=V_{\xi,T}$ be the corresponding weighted
operator. To see a connection of weighted operators with Koopman
representations of compact group extensions consider a compact
(metric) Abelian group $G$ and a cocycle $\va:X\to G$. Then
$U_{T_{\va}}$ (see the glossary) acts on $L^2(X\times G,\mu\ot
\la_G)$. But $$ L^2(X\times G,\mu\ot
\la_G)=\bigoplus_{\chi\in\widehat{G}}L_\chi,\;\;\mbox{where}\;
L_\chi=L^2(X,\mu)\ot\chi,$$ where $L_\chi$ is a
$U_{T_{\va}}$-invariant (clearly, closed) subspace. Moreover, the
map $f\ot\chi\mapsto f$ settles a unitary isomorphism of
$U_{T_{\va}}|_{L_\chi}$ with the operator $V_{\chi\circ \va,T}$.
Therefore, spectral analysis of such Koopman representations
reduces to the spectral analysis of weighted operators
$V_{\chi\circ\va,T}$ for all $\chi\in\widehat{G}$.

\subsection{Maximal spectral type of weighted operators over
rotations}\label{wagowe} The spectral analysis of weighted
operators $V_{\xi,T}$ is especially well developed in case of
rotations, i.e.\ when  we additionally assume that $T$ is an
ergodic rotation on a compact monothetic group $X$: $Tx=x+x_0$,
where $x_0$ is a topologically cyclic element of $X$ (and $\mu$
will stand for Haar measure $\lambda_X$ of $X$). In this case, Helson's
analysis \cite{He} applies (see  also \cite{Gr}, \cite{Iw-Le-Ru},
\cite{Le}, \cite{Qu}) leading to the following
conclusions:\begin{itemize} \item The maximal spectral type
$\sigma_{V_{\xi,T}}$ is either discrete or continuous. \item When
$\sigma_{V_{\xi,T}}$ is continuous it is either singular or
Lebesgue.
\item The spectral multiplicity  of $V_{\xi,T}$ is uniform.
\end{itemize}

We now pass to a description of some results in case when
$Tx=x+\alpha$ is an irrational rotation on the additive circle
$X=[0,1)$. It was already noticed in the original paper by Anzai
\cite{An} that when $\xi:X\to\T$ is an affine cocycle
($\xi(x)=\exp(nx+c)$, $0\neq n\in\Z$) then $V_{\xi,T}$ has a
Lebesgue spectrum. It was then considered by several authors
(originated by  \cite{Ku}, see also \cite{Choe}, \cite{Co-Fo-Si})
to which extent this property is stable when we perturb our
cocycle. Since the topological degree of affine cocyles  is
different from zero, when perturbing them we consider smooth
perturbations by cocycles of degree zero.
\begin{Th}[\cite{Iw-Le-Ru}]\label{we7} Assume that $Tx=x+\alpha$ is an
irrational rotation. If $\xi:[0,1)\to\T$ is of non-zero degree,
absolutely continuous, with the derivative of bounded variation
then $V_{\xi,T}$ has a Lebesgue spectrum.
\end{Th}
In the same paper, it is noticed that if we drop the assumption on
the derivative then the maximal spectral type of $V_{\xi,T}$ is a
Rajchman measure (i.e.\ its Fourier transform vanishes at
infinity). It is still an open question, whether one can find
$\xi$ absolutely continous with non-zero degree and such that
$V_{\xi,T}$ has singular spectrum. ,,Below" absolute continuity,
topological properties of the cocycle seem to stop playing any
role -- in \cite{Iw-Le-Ru} a continuous, degree~1 cocycle $\xi$ of
bounded variation is constructed such that
$\xi(x)=\eta(x)/\eta(Tx)$ for a measurable $\eta:[0,1)\to\T$ (that
is $\xi$ is a {\em coboundary}) and therefore $V_{\xi,T}$ has
purely discrete spectrum (it is isomorphic to $U_T$). For other
results about Lebesgue spectrum for Anzai skew products see also
\cite{Choe}, \cite{Fraczek}, \cite{Iwanik} (in \cite{Fraczek}
$\Z^d$-actions of rotations and so called winding numbers instead
of topological degree are considered). For recent generalizations, see \cite{Ce-Ti}, \cite{Tiedra}.

When the cocycle is still smooth but its degree is zero the
situation drastically changes. Given an absolutely  continuous
function $f:[0,1)\to\R$ M.~Herman \cite{Herman}, using the
Denjoy-Koksma inequality (see e.g.\ \cite{Ku-Ne}), showed that
$f_0^{(q_n)}\to0$ uniformly (here $f_0=f-\int_0^1f\,d\la_{[0,1)}$
and $(q_n)$ stands for the sequence of denominators of $\alpha$).
It follows that $T_{e^{2\pi if}}$ is rigid and hence has a
singular spectrum. B. Fayad  \cite{Fa1} shows that this result is
no longer true if one dimensional rotation is replaced by a
multi-dimensional rotation (his counterexample is in the analytic
class).  See also \cite{Le-Ma} for the singularity of spectrum for
functions $f$ whose Fourier transform satisfies
$\mbox{o}(\frac1{|n|})$ condition or to \cite{Iw-Le-Ma}, where it
is shown that sufficiently small variation implies singularity of
the spectrum.

A natural class of weighted operators arises when we consider
Koopman operators of rotations on nil-manifolds. We only look at
the particular example of such a rotation on a quotient of the
Heisenberg group $(\R^3,\ast)$ (a general spectral theory of
nil-actions was mainly developed by W. Parry \cite{Parry1} --
these actions have countable Lebesgue spectrum in the
orthocomplement of the subspace of eigenfunctions) that is take
the nil-manifold $\R^3/\Z^3$ on which we define the
nil-rotation $
S((x,y,z)\ast\Z^3)=(\alpha,\beta,0)\ast(x,y,z)\ast\Z^3$, where
$\alpha,\beta$ and~1 are rationally independent. It can be shown
that $S$ is isomorphic to the skew product  defined on
$[0,1)^2\times\T$ by
$$
T_\varphi: (x,y,z)\mapsto(x+\alpha,y+\beta,z\cdot e^{2\pi i
\va(x,y)})=(x+\alpha,y+\beta,z+\alpha
y)\ast\Z^3,$$ where $\va(x,y)=\alpha
y-\psi(x+\alpha,y+\beta)+\psi(x,y)$ with $\psi(x,y)=x[y]$. Since
nil-cocycles can be considered as a certain analog of affine
cocycles for one-dimensional rotations, it would be nice to
explain to what kind of perturbations the Lebesgue spectrum
property is stable.

Yet another interesting problem which is related to the spectral
theory of  extensions given by so called {\em Rokhlin cocycles}
(see Section~\ref{czesc8}) arises, when given $f:[0,1)\to\R$, we
want to describe spectrally the one-parameter set of weighted
operators $W_c:=V_{e^{2\pi icf},T}$; here $T$ is a fixed
irrational rotation by~$\alpha$. As proved by quite sophisticated
arguments in \cite{Iw-Le-Ma}, if we take $f(x)=x$ then for all
non-integer $c\in\R$ the spectrum of $W_c$ is continuous and
singular (see also \cite{Gr} and \cite{Me} where some special
$\alpha$'s are considered). It has been open for some time if at
all one can find $f:[0,1)\to\R$ such that for each $c\neq0$, the
operator $W_c$ has a Lebesgue spectrum. The positive answer is
given in \cite{Wy}: for example if $f(x)=x^{-(2+\vep)}$ ($\vep>0$)
and $\alpha$ has bounded partial quotients then $W_c$ has a
Lebesgue spectrum for all $c\neq0$. All functions with such a
property considered in \cite{Wy} are non-integrable. It would be
interesting to find an integrable $f$ with the above property.

We refer to \cite{Go} and the references therein for further
results especially for transformations of the form $(x,y)\mapsto
(x+\alpha,1_{[0,\beta)}(x)+y)$ on $[0,1)\times\Z/2\Z$.  Recall
however that earlier Katok and Stepin  \cite{Ka-St} used this kind
of transformations to give a first counterexample to the
Kolmogorov group property (see the glossary) for the spectrum.

\subsection{The multiplicity problem for weighted operators over rotations}
In case of perturbations of affine cocycles, this problem was
already raised by Kushnirenko  \cite{Ku}. Some significant results
were obtained by M. Guenais. Before we state her results let us
recall a useful criterion to find an upper bound for the
multiplicity: {\em If there exist $c>0$ and a sequence
$(F_n)_{n\geq1}$ of cyclic subspaces of $H$ such that for each
$y\in H$, $\|y\|=1$ we have
$\liminf_{n\to\infty}\|proj_{F_n}y\|^2\geq c$, then $esssup
(M_U)\leq 1/c$} which follows from a well-known lemma of Chacon
(\cite{Chacon}, \cite{Co-Fo-Si}, \cite{Ki}, \cite{Le}). Using this
and a technique which is close to the idea of local rank one (see
\cite{Fe1}, \cite{Ki}) M. Guenais  \cite{Gu1} proved a series of
results on multiplicity which we now list.

\begin{Th} Assume that $Tx=x+\alpha$ and let $\xi:[0,1)\to \T$ be a
cocycle.

(i) If $\xi(x)=e^{2\pi icx}$ then $M_{V_{\xi,T}}$ is bounded by
$|c|+1$.

(ii) If $\xi$ is absolutely continuous and $\xi$ is of topological
degree zero, then $V_{\xi,T}$ has a simple spectrum.

(iii) if $\xi$ is of bounded variation, then
$M_{V_{\xi,T}}\leq\max(2,2\pi Var(\xi)/3)$.
\end{Th}

\subsection{Remarks on the Banach problem}\label{s:Banach} We already mentioned in
Introduction the Banach problem in ergodic theory, which is simply
the question whether there exists a Koopman representation for
$\A=\Z$ with simple Lebesgue spectrum. In fact  no example of a
Koopman representation with Lebesgue spectrum of finite
multiplicity is known. Helson and Parry  \cite{He-Pa} asked for the
validity of a still weaker version: {\em Can one construct $T$
such that $U_T$ has a Lebesgue component in its spectrum whose
multiplicity is finite}? Quite surprisingly in \cite{He-Pa} they
give a general construction yielding for each ergodic $T$  a
cocycle $\va:X\to\Z/2\Z$ such that the unitary operator
$U_{T_{\va}}$ has a Lebesgue spectrum in the orthocomplement of
functions depending only on the $X$-coordinate. Then  Mathew and
Nadkarni \cite{Ma-Na}  gave examples of cocycles over so called
dyadic adding machine for which the multiplicity of the Lebesgue
component was equal to~2. In \cite{Le1} this was generalized to so
called {\em Toeplitz $\Z/2\Z$-extensions} of adding machines: for
each even number $k$ we can find a two-point extension of an
adding machine so that the multiplicity of the Lebesgue component
is~$k$. Moreover, it was shown that Mathew and Nadkarni's
constructions from \cite{Ma-Na} in fact are close to systems
arising from number theory (like the famous Rudin-Shapiro sequence,
e.g.\ \cite{Qu}), relating the result about multiplicity of the
Lebesgue component to results by Kamae \cite{Kamae} and Queffelec
\cite{Qu}. Independently of \cite{Le1}, Ageev  \cite{Ag1} showed
that one can construct 2-point extensions of the Chacon
transformation realizing (by taking powers of the extension) each
even number as the multiplicity of the Lebesgue component.
Contrary to all previous examples, Ageev's constructions are
weakly mixing.

Still an open question remains whether for $\A=\Z$ one can find a
Koopman representation with the Lebesgue component of
multiplicity~1 (or even odd).

In \cite{Gu2}, M. Guenais studies the problem of Lebesgue spectrum
in the classical case of Morse sequences (see \cite{Ke} as well as
\cite{Kw}, where the problem of spectral classification in this
class is studied). All dynamical sytems arising from Morse
sequences have simple spectra \cite{Kw}.  It follows that if a
Lebesgue component appears in a Morse dynamical system, it has
multiplicity one. Guenais \cite{Gu2} using a Riesz product
technique relates the Lebesgue spectrum problem with the still
open problem (a variation of the classical Littlewood problem) of whether a construction of so called $L^1$-ultraflat trigonometric
polynomials with coefficients $\pm1$ is possible (in the very recent preprint \cite{Ba-Bo}, the Littlewood problem of existence of uniformly flat trigonometric polynomials has been solved, but it is unclear whether it yields the ultraflatness condition). However, it is
proved in \cite{Gu2} that such a construction can be carried out
on some compact Abelian groups and it leads, for an Abelian
countable torsion group $\A$, to a construction of an ergodic
action of $\A$ with simple spectrum and a Haar component in its
spectrum.

In
\cite{Pr}, A.\ Prikhodko published a construction of a rank one flow (see Section~\ref{czesc9}) having Lebesgue spectrum. As rank one
implies simple spectrum, the result yields solution of Banach problem for $\A=\R$. To carry out the construction, Prikhodko proved the following $L^1$- ultraflat version of the Littlewood conjecture: {\em For all $0<a<b$
and $n \geq1$, there are polynomials $P_n(t) = \sum_{j=0}^{n-1} e^{2\pi iw_j^{(n)}t}$ for some real numbers $w_j^{(n)}$,
so that $\|P_n\|_{L^1([a,b])}/\sqrt n\to1$ when $n\to\infty$.}
It seems however that some of the arguments in the paper are written too briefly and no further clarifying presentation of methods/results/ideas from \cite{Pr} has appeared so far.\footnote{It would also be extremely nice to explain the status of \cite{Ab} by H.\ El Abdalaoui, first posted on arXiv in 2015, which states the solution of the original Banach problem (i.e.\ in the conservative infinite measure-preserving category).}

\subsection{Lifting mixing properties}\label{liftingmix}
We now give one  more example of interactions between spectral
theory and joinings (see Introduction) that gives rise to a quick
proof of the fact that $r$-fold mixing property of $T$ ($r\geq2$)
lifts to a weakly mixing compact group extension $T_\varphi$ (the
original proof of this fact is due to D. Rudolph \cite{Ru3}).
Indeed, to prove $r$-fold mixing for a mixing(=2-mixing)
transformation $S$ (acting on $\ycn$) one has to prove that each
weak limit of off-diagonal self-joinings (given by powers of $S$,
see \cite{Thierry}) of order $r$ is simply the product measure
$\nu^{\ot r}$. We need also a Furstenberg's lemma (\cite{Fu2})
about relative unique ergodicity (RUE) of compact group extensions
$T_\varphi$: {\em If $\mu\ot \la_G$ is an ergodic measure for
$T_\varphi$ then it is the only (ergodic) invariant measure for
$T_\varphi$ whose projection on the first coordinate is $\mu$.}
Now, the result about lifting $r$-fold mixing to compact group
extensions follows directly from the fact that whenever $T_{\va}$
is weakly mixing, $(\mu\ot \la_G)^{\ot r}$ is an ergodic measure
(this approach was shown to the second author by A.\ del Junco). In particular, if
$T$ is mixing and $T_\varphi$ is weakly mixing then for each
$\chi\in\widehat{G}\setminus\{1\}$, the maximal spectral type of
$V_{\chi\circ\varphi,T}$ is Rajchman.

See Section~\ref{czesc8} for a generalization of the lifting
result to Rokhlin cocycle extensions.

\section{The multiplicity function}\label{czesc6}
In this section only $\A=\Z$ is considered. For other groups, even
for $\R$, much less is known. Clearly, given an automorphism $T$, by inducing its Koopman $\Z$-representation, we obtain a one-parameter group $(V_t)_{t\in\R}$ of unitary operators, which has precisely the same properties as the original one, except that we added the eigenvalues $n\in\R$. Moreover, classically, the induced Koopman representation is given by the suspension of $T$, i.e.\ by the special flow $T^f$ (see the glossary), where $f=1$, whence it is also Koopman but is never weakly mixing. See  \cite{Da-Le100} and \cite{Da-So}, where some of the results below  proved for $\A=\Z$ have been extended to (weakly mixing) flows. See also  the case of so called
{\em product} $\Z^d$-{\em actions} \cite{Fi} and \cite{So} for general countable Abelian group actions.

Contrary to the
case of maximal spectral type, it is rather commonly believed that
there are no restrictions for the set of essential values of
Koopman representations. In fact, if we drop the assumption that we consider the finite measure-preserving case and let ourselves consider $\mu$ $\sigma$-finite and infinite, Danilenko and Ryzhikov \cite{Da-Ry}, \cite{Da-Ry1} proved that all subsets of $\{1,2,\ldots\}\cup\{\infty\}$ are Koopman realizable (in the weak mixing and mixing class, respectively).

\subsection{Cocycle approach}
We will only concentrate on some results of the last four decades. In 1983, refining an earlier idea of Oseledets from the 1960th, E.A.\ Robinson \cite{Ro1} proved that for each
$n\geq1$ there exists an ergodic transformation whose maximal
spectral multiplicity is~$n$.  Another important result was proved
in \cite{Ro3} (see also \cite{Ka0}), where it is shown that given
a finite set $M\subset\N$ containing~1 and closed under the least
common multiple one can find (even a weakly mixing) $T$ so that
the set of essential values of the multiplicity function equals
$M$. This result was then extended in \cite{Go-Kw-Le-Li} to
infinite sets and finally in \cite{Kwjunior-Le} (see also
\cite{Ag2}) to all subsets $M\subset\N$ containing~1. In fact, as
we have already noticed in the previous section the spectral
theory for compact Abelian group extensions is reduced to a study of
weighted operators and then to comparing maximal spectral types for
such operators.  This leads to sets of the form
$$
M(G,v,H)=\{\sharp(\{\chi\circ v^i:\:i\in\Z\}\cap
\,anih(H)):\:\chi\in\,anih(H)\}
$$
($H\subset G$ is a closed subgroup and $v$ is a continuous group
automorphism of $G$).
 Then an algebraic lemma has been proved in \cite{Kwjunior-Le}
saying that each set $M$ containing~1 is of the form $M(G,v,H)$
and the techniques to construct ``good" cocycles and a passage to
``natural factors" yielded the following: {\em For each
$M\subset\{1,2,\ldots\}\cup\{\infty\}$ containing~1 there exists
an ergodic automorphim such that the set of essential values for
its Koopman representation equals~$M$.} See also \cite{Ro2} for
the case of non-Abelian group extensions.

A similar in spirit approach (that means, a passage to a family of
factors) is present in the paper of Ageev \cite{Ag3} in which
he first applies Katok's analysis (see \cite{Ka0}, \cite{Ka01})
for spectral multiplicities of the Koopman representation
associated with Cartesian products $T^{\times k}$ for a generic
transformation $T$. In a natural way this approach leads to study
multiplicities of tensor products of unitary operators. Roughly,
the multiplicity is computed as the number of atoms (counted
modulo obvious symmetries) for conditional measures (see
\cite{Ka0}) of a product measure over its convolution. Ageev
\cite{Ag3} proved  that for a typical automorphism $T$ the set of
the values of the multiplicity function for $U_{T^{\times k}}$
equals $\{k,k(k-1),\ldots, k!\}$ and then he just passes to
``natural" factors for the Cartesian products by taking sets
invariant under a fixed subgroup of permutations of coordinates.
In particular, he obtains all sets of the form $\{2,3,\ldots,n\}$
on $L^2_0$. He also shows that such sets of multiplicities are
realizable in the category of mixing transformations. See also Ryzhikov \cite{Ry2009} for a realization of sets of the form $\{k,l,kl\}$, $\{k,l,m,kl,km,lm,klm\}$, etc.

A further progress was done in 2009-2012, when first Katok and Lema\'nczyk \cite{Ka-Le} proved that each finite subset $M\subset\{1,2,\dots\}\cup\{\infty\}$ containing~2 can be realized as the set of essential values of an ergodic automorphism which was then, by overcoming some algebraic difficulties, extended by Danilenko \cite{Da0}, \cite{Da00} (in the mixing category) to all subsets containing~2.

\subsection{Multiplicity for Gaussian and Poissonian automorphisms} We refer the reader to Section~\ref{czesc10} for the definition and basic properties of Gaussian and Poissonian automorphisms. Recall that given a Poissonian automorphism, there is a Gaussian automorphism spectrally isomorphic to it (whether the converse holds, is unknown). In Gaussian case, the classical Girsanov's theorem from the 1950th asserts that the maximal spectral multiplicity in this case is either one or infinity (with a possibility that $\infty$ is not an essential value), see \cite{Ku-Pa} for an elegant proof of this theorem. What is the family of subsets appearing as sets of essential values of the multiplicity functions of Koopman operators given by Gaussian (and also Poissonian)  automorphisms was studied by Danilenko and Ryzhikov in \cite{Da-Ry1}. They prove the following remarkable results:
\begin{itemize} \item
this family contains all  multiplicative sub-semigroups of $\N$; \item this family contains other sets than multiplicative sub-semigroups of $\N$.
\end{itemize}
The latter shows that Proposition 6.4.4 (multiplicative nature of $M_T$ in the Gaussian case) claimed in the book \cite{Ka-Th} (and also in \cite{Ro}) is not true.

In the unpublished preprint \cite{RyArx}, Ryzhikov shows that all subsets containing $\infty$ are Gaussian ``realizable'' (even in the mixing category).

\subsection{Rokhlin's uniform multiplicity problem}
The Rokhlin multiplicity problem (see the recent book by Anosov
\cite{Anosov}) was, given $n\geq2$, to construct an ergodic
transformation with uniform multiplicity~$n$ on $L^2_0$. A
solution for $n=2$ was independently given by Ageev \cite{Ag4} and
Ryzhikov \cite{Ry3} (see also \cite{Anosov} and \cite{Go}) and in
fact it consists in showing that for some $T$ (actually, any $T$
with simple spectrum for which $\frac12(Id+U_T)$ is in the weak
operator closure of the powers of $U_T$ will do) the multiplicity
of $T\times T$ is uniform and equal to~2 (see also
Section~\ref{czesc13}). In  \cite{Ti0} (\cite{Ry-Tr}), the case $n=2$ is solved in case of mixing automorphisms (flows).

In \cite{Ag5}, Ageev proposed a new approach which consists in
considering actions of ,,slightly non-Abelian" groups and showing
that for a ``typical" action of such a group a fixed ``direction"
automorphism has a uniform multiplicity. Shortly after publication
of \cite{Ag5}, Danilenko  \cite{Da}, following Ageev's approach,
considerably simplified the original proof. We will present
Danilenko's arguments.

 Fix $n\geq1$. Denote
$\ov{e}_i=(0,\ldots,1,\ldots,0)\in \Z^n$, $i=1,\ldots,n$. We
define an automorphism  $L$ of $\Z^n$ setting $L(\ov{e}_i)=
\ov{e}_{i+1}$, $i=1,\ldots,n-1$ and $L(\ov{e}_n)=\ov{e}_1$. Using
$L$ we define a semi-direct product  $G:=\Z^n \rtimes \Z$ defining
multiplication as $(u,k)\cdot(w,l)=(u+L^kw,k+l)$. Put $e_0=(0,1)$,
$e_i=(\ov{e}_i,0)$, $i=1,\ldots,n$ (and $Le_i=(L\ov{e}_i,0)$).
Moreover, denote $e_{n+1}=e_0^n=(0,n)$. Notice that $e_0\cdot
e_i\cdot e_0^{-1}=Le_i$ for $i=1,\ldots,n$ ($L(e_{n+1})=e_{n+1}$).

\begin{Th}[Ageev, Danilenko]\label{agiejew2}
For every unitary representation  $\cu$ of $G$ in a separable
Hilbert space $H$, for which $U_{e_1-L^re_1}$ has no non-trivial
fixed points for $1\leq r<n$, the essential values of the
multiplicity function for  $U_{e_{n+1}}$ are contained in the set
of multiples of~$n$.  If, in addition, $U_{e_0}$ has a simple
spectrum, then $U_{e_{n+1}}$ has uniform multiplicity~$n$.
\end{Th}

It is then a certain work to show that the assumption of the
second part of the theorem  is satisfied for a typical action of
the group $G$. Using a special $(C,F)$-construction with all the
cut-and-stack parameters explicit,  Danilenko  \cite{Da} was also
able to show that each set of the form $k\cdot M$, where $k\geq1$
and $M$ is an arbitrary subset of natural numbers containing~1, is
realizable as the set of essential values of a Koopman
representation.

Tikhonov \cite{Ti0} proved the existence of a mixing automorphism of uniform multiplicity $n$ on $L^2_0$ for all $n\geq1$.

\section{Rokhlin cocycles}\label{czesc8}
We consider now a certain class of extensions which should be
viewed as a generalization of the concept of compact group
extensions. We will focus on $\Z$-actions only.

Assume that $T$ is an ergodic automorphism of $\xbm$. Let $G$ be a
l.c.s.c. Abelian  group. Assume that this group acts on $\ycn$,
that is we have a $G$-action $\cs=(S_g)_{g\in G}$ on $\ycn$. Let
$\va:X\to G$ be a cocycle. We then define an automorphism
$T_{\va,\cs}$ of the space $(X\times Y,\cb\ot\cc,\mu\ot\nu)$ by
$$
T_{\va,\cs}(x,y)=(Tx,S_{\va(x)}(y)).
$$
Such an extension is called a {\em Rokhlin cocycle extension} (the
map $x\mapsto S_{\va(x)}$ is called a {\em Rokhlin cocycle}). Such
an operation generalizes the case of compact group extensions;
indeed, when  $G$ is  compact the action of $G$ on itself by
rotations preserves Haar measure. (It is quite surprising, that
when only we admit $G$ non-Abelian, then, as shown in
\cite{Da-Le}, {\bf each} ergodic extension of $T$ has a form of a
Rokhlin cocycle extension.) Ergodic and spectral properties of
such extensions are examined in several papers:
\cite{Gl1}, \cite{Gl-We}, \cite{Le-Le}, \cite{Le-Me-Na},
\cite{Le-Pa1}, \cite{Le-Pa2}, \cite{Ro4}, \cite{Ru4}. Since in
these papers rather joining aspects are studied (among other
things in \cite{Le-Le} Furstenberg's RUE lemma is generalized to
this new context), we will mention here only few results, mainly
spectral, following \cite{Le-Le} and \cite{Le-Pa2}. We will
constantly assume that $G$ is non-compact. As $\va:X\to G$ is then
a cocycle with values in a non-compact group, the theory of such
cocycles is much more complicated (see e.g.\ \cite{Sch}), and in
fact the theory of Rokhlin cocycle extensions leads to interesting
interactions between classical ergodic theory, the theory of
cocycles  and the theory of non-singular actions arising from
cocycles taking values in non-compact groups -- especially, the Mackey
action associated to $\va$ plays a crucial role here (see the
problem of invariant measures for $T_{\va,\cs}$ in \cite{Le-Pa1}
and \cite{Da-Le}); see also monographs \cite{Aa}, \cite{Ka0},
\cite{Ka02}, \cite{Sch}. Especially, two Borel subgroups of
$\widehat{G}$ are important here:
$$
\svf=\{\chi\in\widehat{G}:\:\chi\circ\va=c\cdot\xi/\xi\circ
T\;\mbox{for a measurable}\;\xi:X\to\tor\;\mbox{and}\;c\in\tor\}.
$$
and its subgroup $\Lambda_\varphi$ given by $c=1$. $\lf$ turns out
to be the group of $L^\infty$-eigenvalues of the Mackey action (of
$G$) associated to the cocycle $\va$. This action is the quotient
action of the natural action of $G$ (by translations on the second
coordinate) on the space of ergodic components of the skew product
$T_{\va}$ -- the Mackey action is (in general) not
measure-preserving, it is however non-singular.  We refer the
reader to \cite{Aa-Na}, \cite{Ho-Me-Pa} and \cite{Na} for other
properties of those subgroups.

\begin{Th}[\cite{Le-Pa1},\cite{Le-Pa2}]\label{p1}
(i) $\stfs|_{L^2(X\times Y,\mu\otimes\nu)\ominus
L^2(X,\mu)}=\int_{\hat{G}}\sigma_{V_{\chi\circ\va,T}}\,d\sigma_{\cs}$.

(ii) $\tfs$ is ergodic if and only if $T$ is ergodic and
$\si_{\cs}(\lf)=0$.

(iii) $\tfs$ is weakly mixing if and only if $T$ is weakly mixing
and $\cs$ has no eigenvalues in $\svf$.

(iv) if $T$ is mixing, $\cs$ is mildly mixing, $\va$ is recurrent
and not cohomologous to a cocycle with values in a compact
subgroup of $G$ then $T_{\va,\cs}$ remains mixing.

(v) If $T$ is $r$-fold mixing, $\va$ is recurrent and
$T_{\va,\cs}$ is mildly mixing then $T_{\va,\cs}$ is also $r$-fold
mixing.

(vi) If $T$ and $R$ are disjoint, the cocycle $\va$ is ergodic and
$\cs$ is mildly mixing then $T_{\va,\cs}$ remains disjoint from
$R$.
\end{Th}

Let us recall (\cite{Fu-We}, \cite{Sch-Wa}) that an $\A$-action
$\cs=(S_a)_{a\in\A}$ is mildly mixing (see the glossary) if and
only if the  $\A$-action $(S_a\times\tau_a)_{a\in\A}$ remains
ergodic for every properly ergodic non-singular  $\A$-action
$\tau=(\tau_a)_{a\in\A}$.

Coming back to Smorodinsky-Thouvenot's result about factors of
ergodic self-joinings of a Bernoulli automorphism we would like to
emphasize here that the disjointness result~(vi) above was used in
\cite{Le-Pa1} to give an example of an automorphism which is
disjoint from all weakly mixing transformations but which has an
ergodic self-joining whose associated automorphism  has a
non-trivial weakly mixing factor. In a sense this is opposed to
Smorodinsky-Thouvenot's result as here from self-joinings we
produced a ``more complicated"  system (namely the weakly mixing
factor) than the original system.

It would be interesting to develop the theory of spectral
multiplicity for Rokhlin cocycle extensions as it was done in the case
of compact group extensions. However a difficulty is that in the
compact group extension case we deal with a countable direct sum
of representations of the form $V_{\chi\circ\va,T}$ while in the
non-compact case we have to consider an integral of such
representations.

\section{Rank-1 and related systems}\label{czesc9}
 Although properties like mixing, weak
(and mild) mixing as well as ergodicity, are clearly spectral
properties, they have ``good" measure-theoretic formulations
(expressed by a certain behaviour on sets). Simple spectrum
property is another example of a spectral property, and it was a
popular question in the 1980s whether simple spectrum property of
a Koopman representation can be expressed in a more
``measure-theoretic" way. We now recall rank-1 concept which can
be seen as a notion close to Katok's and Stepin's theory of cyclic
approximation \cite{Ka-St} (see also \cite{Co-Fo-Si}).

Assume that $T$ is an automorphism of a standard probability Borel
space $\xbm$. $T$ is said to have the {\em rank one} property if there
exists an increasing sequence of Rokhlin towers tending to the
partition into points (a {\em Rokhlin tower} is a family
$\{F_n,TF_n,\ldots,T^{h_n-1}F_n\}$ of pairwise disjoint sets, while
``tending to the partition into points" means that we can
approximate every set in $\cb$ by unions of levels of towers in
the sequence). Hence, basically, rank one automorphism is given
by two sequences of parameters: $r_n$, $n\geq1$, which is the number of subcolumns
on which we divide the $n$th tower given by $F_n$, and $s_{n,j}$, $n\geq1$, $j=0,1,\ldots,r_n-1$, the sequence of spacers put over consecutive subcolumns. A ``typical'' automorphism of a standard probability Borel space  has the rank one property. The ``typicality'' of rank one is still true in the Alpern-Tikhonov topology we mentioned in Section~\ref{czesc4} by \cite{Bash1}.

Baxter  \cite{Ba} showed that the maximal spectral
type of such a $T$ is realized by a characteristic function. Since
the cyclic space generated by the characteristic function of the
base contains characteristic functions of all levels of the tower,
by the definition of rank one, the increasing sequence of cyclic
spaces tends to the whole $L^2$-space, therefore rank one property
implies simplicity of the spectrum for the Koopman representation.
It was a question for some time whether rank-1 is just a
characterization of simplicity of the spectrum, disproved by del
Junco \cite{Ju2}. We refer the reader to \cite{Fe3} as a good
source for basic properties of rank-1 transformations.

Similarly to the rank one property, one can define {\em finite
rank} automorphisms (simply by requiring that an approximation is
given by a sequence of a fixed number of towers) -- see e.g.
\cite{Or-Ru-We}, or even, a more general property, namely the {\em
local rank one} property can be defined, just by requiring that
the approximating sequence of single towers fills  up a fixed
fraction of the space (see \cite{Fe1}, \cite{Ki}). Local rank one
property implies finite multiplicity
\cite{Ki} and the maximal spectral multplicity is always bounded by rank. Mentzen \cite{Men} showed that for each $n\geq1$ one
can construct an automorphism with simple spectrum and having
rank~$n$ and later Kwiatkowski and Lacroix \cite{Kw-La} showed that for each pair $(m,r)$ with $m\leq r$, one can construct a rank r automorphism whose maximal spectral multplicity is~$m$.  In \cite{Le-Si} there is an example of a simple spectrum
automorphism which is not of local rank one. Ferenczi \cite{Fe2}
deals with the notion of funny rank one (approximating towers are
Rokhlin towers with ``holes") - the concept that has been introduced by Thouvenot. Funny rank one also implies
simplicity of the spectrum. An example is given in \cite{Fe2}
which is even not loosely Bernoulli (see Section~\ref{czesc11}, we
recall that local rank one property implies loose Bernoullicity
\cite{Fe1}).

The notion of AT (see the glossary) has been introduced by Connes
and Woods~\cite{Co-Wo}. They proved that AT property implies zero
entropy. They also proved that funny rank one automorphisms are
AT. In \cite{Do-Qu} it is proved that the system induced by the
classical Morse-Thue system is AT (it is an open question by S.
Ferenczi whether this system has funny rank one property). A
question by Dooley and Quas is whether AT implies funny rank one
property. AT property implies ``simplicity of the spectrum in
$L^1$" which we already considered in Introduction (a ``generic"
proof of this fact is due to J.-P. Thouvenot).

A persistent question was formulated in the 1980s whether rank one
itself is a spectral property. In \cite{Fe-Le}, the authors
maintained that this is not the case, based on an unpublished
preprint of the first named author of~\cite{Fe-Le} in which there
was a construction of a Gaussian-Kronecker  automorphism (see
Section~\ref{czesc10}) having rank-1 property. This latter
construction turned out to be false. In fact, de la Rue
\cite{delaRue3} proved that no Gaussian automorphism can be of
local rank one. Therefore the question whether: {\em Rank one is a
spectral property} remains one of interesting open questions in
that theory. Downarowicz and Kwiatkowski \cite{Do-Kw1} proved that
rank-1 is a spectral property in the class of systems generated by
generalized Morse sequences.

One of the most beautiful theorems about rank-1 automorphisms is
the following result of J. King \cite{Ki5} (for a different proof
see \cite{Ry7}).

\begin{Th}[WCT] If $T$ is of rank one then
for each element $S$ of the centralizer $C(T)$ of $T$ there exists
a sequence $(n_k)$ such that $U^{n_k}_T\to U_S$ strongly.
\end{Th}

A conjecture of J. King is that in fact for rank-1 automorphisms
each indecomposable Markov operator $J=J_\rho$ ($\rho\in
J_2^e(T)$) is a weak limit of powers of $U_T$ (see \cite{Ki6},
also \cite{Ry7}). To which extent the WCT remains true for actions
of other groups is not clear. In \cite{Ze} the WCT is proved in
case of rank one flows, however the main argument seems to be
based on the fact that a rank one flow has a non-zero time
automorphism $T_{t_0}$ which is of rank one, which is not true.
After the proof of the WCT by Ryzhikov in \cite{Ry7} there is a
remark that the rank one flow version of the theorem can be proved
by a word for word repetition of the arguments. He also proves
that if the flow $(T_t)_{t\in\R}$ is mixing, then $T_1$ does not
have finite rank. On the other hand, for $\A=\Z^2$, Downarowicz
and Kwiatkowski \cite{Do-Kw2} gave a counterexample to
the WCT. But see also \cite{Ja-Ru-Ry}.

Even though it looks as if rank one construction is not
complicated, mixing in this class is possible; historically the
first mixing constructions were given by D. Ornstein \cite{Or} in
1970, using probability type arguments for a choice of spacers.
Once mixing was shown, the question arose whether absolutely
continuous spectrum is also possible, as this would give
automatically the positive answer to the Banach problem. However
Bourgain \cite{Bo}, relating spectral measures of rank one
automorphisms with some classical constructions of Riesz product
measures, proved that a certain subclass of Ornstein's class
consists of automorphisms with singular spectrum (see also
\cite{Ab1} and \cite{Ab-Pa-Pr}). Since in Ornstein's class spacers
are chosen in a certain ``non-constructive" way, quite a lot of
attention was devoted to the rank one automorphism defined by
cutting a tower at the $n$-th step into $r_n=n$ subcolumns of
equal ``width" and placing $i$ spacers over the $i$-th subcolumn.
The mixing property conjectured by M. Smorodinsky,  was proved by
Adams \cite{Ad} (in fact Adams proved a general result on mixing
of a class of staircase transformations). Spectral properties of
rank-1 transformations are also studied  in \cite{Kl-Re}, where
the authors proved that whenever $\sum_{n=1}^\infty
r_n^{-2}=+\infty$  then
the spectrum is automatically singular, see also more recent \cite{Cr-Si}.   H. Abdalaoui \cite{Ab1}
gives a criterion for singularity of the spectrum of a rank one
transformation; his proof uses a central limit theorem. It seems
that still the question whether rank one implies singularity of
the spectrum remains the most important question of this theory.

We have already seen in Section~\ref{czesc5} that for a special
class of rank one systems, namely those with discrete spectra
(\cite{Ju1}), we have a nice theory for weighted operators. It
would be extremely interesting to find  a rank one automorphism with
continuous spectrum for which  a substitute of  Helson's analysis
exists.

B. Fayad \cite{Fa} constructs a rank one differentiable flow, as a
special flow over a two-dimensional rotation. In \cite{Fa2} he
gives new constructions of smooth flows with singular spectra
which are mixing (with a new criterion for a Rajchman measure to
be singular). In \cite{Fa3}  a certain smooth change of time for
an irrational flows on the 3-torus is given, so that the
corresponding flow  is partially mixing and has the local rank one
property.

Motivated by Sarnak's conjecture on M\"obius disjointness, see \cite{Ku-LeSurvey}, a certain recent activity was to study spectral disjointness of powers for rank one automorphisms.
Let $\sigma$ be a probability measure on the additive circle $[0,1)$. Given a real number $a>0$, we denote by
$\sigma^a$ the image of $\sigma$ under the map $x\mapsto ax$ mod~1. If $r\geq1$ is an integer, then
by $\sigma_r$, we will denote the measure which is obtained first by taking the image of $\sigma$ under the map $x\mapsto\frac1rx$, i.e.\ the measure $\sigma^{1/r}$, and then repeating this new measure periodically in intervals $[\frac jr,\frac{j+1}r)$. The following holds:\footnote{e.g. in the unpublished notes by H. El Abdalaoui, J.\ Ku\l aga-Przymus, M.\ Lema\'nczyk and T. de la Rue.}
\begin{equation}\label{l:bourgaintype}
\mbox{if $(r,s)=1$ then $\sigma^r\perp\eta^s$ if and only if $\sigma_s\perp \eta_r$.}
\end{equation}

In \cite{Bo1}, Bourgain used Riesz product technique to show that for the class of so called rank one automorphisms with bounded parameters (both $(r_n)$ and $(s_{n,j})$ are bounded and no spacer over the last column) we have $\sigma_r\perp\sigma_s$ for $r\neq s$ prime. In view of~(\ref{l:bourgaintype}) it follows that different prime powers are spectrally disjoint.
In \cite{Ab-Le-Ru}, a much larger class of rank one automorphisms is considered. No boundedness assumption on $(r_n)$ is made but a certain bounded recurrence is required on the sequence of spacers. Spectral disjointness of different powers (for the continuous part of the maximal spectral type) is derived from the existence, in the weak closure of powers, of sufficiently many analytic functions of the Koopman operator $U_T$.

For a spectral disjointness of the continuous part of the maximal spectral type for powers of automorphisms like the substitutional system  given by the Thue-Morse sequence and related (rank two systems), see \cite{Ab-Ka-Le}. Weak closure of powers for Chacon automorphism is described in \cite{Ja-Pr-Ru-Ry}.

\section{Spectral theory of dynamical systems of probabilistic
origin}\label{czesc10} Let us just recall that when
$(Y_n)_{n=-\infty}^\infty$ is a stationary process then its
distribution $\mu$ on $\R^{\Z}$  is invariant under the shift $S$
on $\R^{\Z}$: $S((x_n)_{n\in\Z})=(y_n)_{n\in\Z}$, where
$y_n=x_{n+1},\;n\in\Z$. In this way we obtain an automorphism $S$
defined on $(\R^{\Z},\cb(\R^{\Z}),\mu)$. For each automorphism $T$
we can find $f:X\to\R$ measurable such that the smallest
$\sigma$-algebra making the stationary process $(f\circ
T^n)_{n\in\Z}$ measurable is equal to $\cb$, therefore, for the
purpose of this article,  by a system of probabilistic origin we
will mean $(S,\mu)$ obtained from a stationary infinitely
divisible process (see e.g.\ \cite{Mar}, \cite{Sa}). In
particular, the theory of Gaussian dynamical systems is indeed a
classical part of ergodic theory (e.g.\  \cite{Ne}, \cite{ParryG},
\cite{Ver1}, \cite{Ver2}). If $(X_n)_{n\in\Z}$ is a stationary
real centered Gaussian process and $\sigma$ denotes the {\em
spectral measure of the process}, i.e.\
$\widehat{\sigma}(n)=E(X_n\cdot X_0)$, $n\in\Z$, then by
$S=S_\sigma$ we denote the corresponding Gaussian system on the
shift space (recall also that for each symmetric measure $\sigma$
on $\T$ there is exactly one stationary real centered Gaussian
process whose spectral measure is $\sigma$). Notice that if
$\sigma$ has an atom, then in the cyclic space generated by $X_0$
there exists an eigenfunction $Y$ for $S_\sigma$ -- if now
$S_\sigma$ were ergodic, $|Y|$ would be a constant function which
is not possible by the nature of elements in $\Z(X_0)$. In what
follows we assume that $\sigma$ is continuous.

It follows that $U_{S_\sigma}$ restricted to $\Z(X_0)$ is
spectrally the same as $V=V^\sigma$ acting on $L^2(\T,\sigma)$,
and we obtain that  $(U_{S_\sigma},L^2(\R^{\Z},\mu_\sigma))$ can
be represented as the symmetric Fock space built over
$H=L^2(\T,\sigma)$ and $U_{S_\sigma}=F(V)$ -- see the glossary
($H^{\odot n}$ is called the $n$-{\em th chaos}). In other words
the spectral theory of Gaussian dynamical systems is reduced to
the spectral theory of special tensor products unitary operators.
Classical results (see \cite{Co-Fo-Si}) which can be obtained from
this point of view are for example the following:\\
(A) ergodicity implies weak mixing,\\
(B) the multiplicity function is either~1 or is unbounded,\\
(C) the maximal spectral type of $U_{S_\sigma}$ is equal to
$\exp(\sigma)$, hence Gaussian systems enjoy the Kolmogorov group
property.

However, we can also look at a Gaussian system in a different way,
simply by noticing that the  variables $e^{2\pi if}$ ($f$ is a
real variable), where $f\in\Z(X_0)$ generate
$L^2(\R^{\Z},\mu_\sigma)$. Now, calculating
 the spectral measure of $e^{2\pi if}$ is not difficult
and we obtain easily~(C). Moreover, integrals of type $\int
e^{2\pi if_0}e^{2\pi if_1\circ T^n}e^{2\pi if_2\circ
T^{n+m}}\,d\mu_\sigma$ can also be calculated, whence in
particular we easily obtain Leonov's theorem on the multiple
mixing property of Gaussian systems \cite{Leonov}.

One of the most beautiful parts of the theory of Gaussian systems
concerns ergodic  properties of $S_\sigma$  when $\sigma$ is
concentrated on a thin Borel set. Recall that a closed subset
$K\subset\T$ is said to be a {\em Kronecker set} if each $f\in
C(K)$ is a uniform limit of characters (restricted to $K$). Each
Kronecker set has no rational relations.  Gaussian-Kronecker
automorphisms are, by definition, those Gaussian systems for which
the measure $\sigma$ (always assumed to be continuous) is
concentrated on $K\cup\ov{K}$, $K$ a Kronecker set. The following
theorem has been proved in \cite{Fo-St} (see also
\cite{Co-Fo-Si}).

\begin{Th}[Foia\c{s}-Stratila Theorem]
If $T$ is an ergodic automorphism and $f$ is a real-valued element
of $L^2_0$ such that the spectral measure $\sigma_{f}$ is
concentrated on $K\cup\ov{K}$, where $K$ is a Kronecker set, then
the process $(f\circ T^n)_{n\in\Z}$ is Gaussian.
\end{Th}

This theorem is indeed striking as it gives examples of weakly
mixing automorphisms which are spectrally determined (like
rotations). A relative version of the Foia\c{s}-Stratila Theorem
has been  proved in \cite{Le-Le}.

The Foia\c{s}-Stratila Theorem implies that whenever a spectral
measure $\sigma$ is Kronecker, it has no realization of the form
$\sigma_f$ with $f$  bounded. We will see however in
Section~\ref{czesc13} that for some automorphisms $T$ (having the
SCS property) the maximal spectral type $\sigma_T$ has the
property that $S_{\sigma_T}$ has a simple spectrum.

Gaussian-Kronecker automorphisms are examples of automorphisms
with simple spectra. In fact, whenever $\sigma$ is concentrated on
a set without rational relations, then $S_\sigma$ has a simple
spectrum (see \cite{Co-Fo-Si}). Examples of mixing automorphisms
with simple spectra are known \cite{Ne}, however it is still
unknown (Thouvenot's question) whether the Foia\c{s}-Stratila
property may hold in the mixing class. F.\ Parreau \cite{Parreau}
using independent Helson sets gave an example of mildly mixing
Gaussian system with the Foia\c{s}-Stratila property.

Joining theory of  a class of Gaussian system, called GAG, is
developed in \cite{Le-Pa-Th}. A Gaussian automorphism $S_\sigma$
with the Gaussian space $H\subset L^2_0(\R^{\Z},\mu_\sigma)$ is
called a GAG if for each ergodic self-joining $\rho\in
J_2^e(S_\sigma)$ and arbitrary $f,g\in H$ the variable $$
(\R^{\Z}\times \R^{\Z},\rho)\ni(x,y)\mapsto f(x)+g(y)$$ is
Gaussian. For GAG systems one can describe the centralizer and
factors, they turn out to be objects close to the probability
structure of the system. One of the crucial observations in
\cite{Le-Pa-Th} was that all Gaussian systems with simple spectrum
are GAG.

It is conjectured (J.P. Thouvenot) that in the class of zero
entropy Gaussian systems the PID property holds true.

For the spectral theory of classical factors of a Gaussian system
see \cite{Le-Sam de Lazaro}; also spectrally they share basic
spectral properties of Gaussian systems. Recall that historically
one of the classical factors namely the $\sigma$-algebra of sets
invariant for the map $$ (\ldots,x_{-1},x_0,x_1,\ldots)\mapsto
(\ldots,-x_{-1},-x_0,-x_1,\ldots)$$ was the first example with
zero entropy and countable Lebesgue spectrum (indeed, we need a
singular measure $\sigma$ such that $\sigma\ast\sigma$ is
equivalent to Lebesgue measure \cite{ParryG}). For factors
obtained as functions of a stationary process see
\cite{Iw-Le-La-Ru}.

T. de la Rue   \cite{delaRue3} proved  that Gaussian systems are
never of local rank-1, however his argument does not apply to
classical factors. We conjecture that Gaussian systems are
disjoint from rank-1 automorphisms (or even from local rank-1
systems).

We now turn the attention to Poissonian systems (see
\cite{Co-Fo-Si} for more details). Assume that $\xbm$ is a
standard Borel space, where $\mu$ is infinite. Without entering too much into details, the new
configuration space $\widetilde{X}$ is taken as the set of all
countable subsets $\{x_i:\:i\geq1\}$ of $X$. Once a set $A\in\cb$,
of finite measure is given one can define a map
$N_A:\widetilde{X}\to\N(\cup\{\infty\})$ just counting the number
of elements belonging to $A$. The  measure-theoretic structure
$(\widetilde{X},\widetilde{\cb},\widetilde{\mu})$ is given so that
the maps $N_A$ become random variables with Poisson distribution
of parameter $\mu(A)$ and such that whenever
$A_1,\ldots,A_k\subset X$ are of finite measure and are pairwise
disjoint then the variables $N_{A_1},\ldots,N_{A_k}$ are
independent.

Assume now that $T$ is an automorphism of $\xbm$. It induces a
natural automorphism on the space
$(\widetilde{X},\widetilde{\cb},\widetilde{\mu})$ defined by
$\widetilde{T}(\{x_i:\:i\geq1\}=\{Tx_i:\:i\geq1\}$. The
automorphism $\widetilde{T}$ is called the {\em Poisson
suspension} of $T$ (see \cite{Co-Fo-Si}).  Such a suspension is
ergodic if and only if no set of positive and finite $\mu$-measure
is $T$-invariant. Moreover, the ergodicity of $\widetilde{T}$ implies
weak mixing. In fact the spectral structure of $U_{\widetilde{T}}$
is very similar to the Gaussian one: namely the first chaos equals
$L^2\xbm$ (we emphasize that this is about the whole $L^2$ and not
only $L^2_0$) on which $U_{\widetilde{T}}$ acts as $U_T$ and the
$L^2(\widetilde{X},\widetilde{\mu})$ together with the action of
$U_{\widetilde{T}}$ has the structure of the symmetric Fock space
$F(L^2\xbm)$ (see the glossary).

We refer  to \cite{Ca-Po-We}, \cite{Ja-We}, \cite{Ros-Zak1},
\cite{Ros-Zak2} for ergodic properties of systems given by
symmetric $\alpha$-stable stationary processes, or more generally
infinitely divisible processes. Again, they share spectral
properties similar to the Gaussian case: ergodicity implies weak
mixing, while mixing implies mixing of all orders.

In \cite{Roy}, E. Roy clarifies the dynamical ``status" of such
systems. He uses Poisson suspension automorphisms and the Maruyama
representation of an infinitely divisible process mixed with basic
properties of automorphisms preserving infinite measure (see
\cite{Aa}) to prove that as a dynamical system, a stationary
infinitely divisible process (without the Gaussian part), is a
factor of the Poisson suspension over the L\'evy measure of this
process. In \cite{Roy2} a theory of ID-joinings is developed
(which should be viewed as an analog of the GAG theory in the
Gaussian class). Parreau and Roy \cite{Pa-Ro} study
Poisson suspensions without non-trivial factors.

Many natural problems still remain open here, for example
(assuming always zero entropy of the dynamical system under
consideration): Are Poisson suspensions disjoint from Gaussian
systems? In \cite{Ja-Ro-Ru} there are examples of Poissonian systems which are disjoint from all Gaussian systems. What is the spectral structure for dynamical systems
generated by symmetric $\alpha$-stable processes? Are such systems
disjoint whenever $\alpha_1\neq\alpha_2$?  Are Poissonian systems
disjoint from local rank one automorphisms (cf.\ \cite{delaRue3})?
In \cite{Ju-Le2} it is proved that Gaussian systems are disjoint
from so called simple systems (see \cite{Ve0}, \cite{Ju-Ru} and
\cite{Thierry}); we will come back to an extension of this result
in Section~\ref{czesc13}. It seems that flows of probabilistic
origin satisfy the Kolmogorov group property for the spectrum. One
can therefore ask how different  are systems satisfying the
Kolmogorov group property from systems for which the convolutions
of the maximal spectral type are pairwise disjoint (see also
Section~\ref{czesc13} and the SCS property).

We also mention here another problem which will be taken up in
Section~\ref{czesc12} -- {\em Is it true that flows of
probabilistic origin are disjoint from smooth flows on surfaces?}

Yet one more (joining) property seems to be characteristic in the
class of systems of probabilistic origin, namely they satisfy so
called ELF property (see \cite{De-Fr-Le-Pa} and de la Rue's
article \cite{Thierry}). Vershik asked whether the ELF property is
spectral -- however the example of a cocycle from \cite{Wy}
together with Theorem~\ref{p1}~(i) yields a certain Rokhlin
extension of a rotation which is ELF and has countable Lebesgue
spectrum in the orthocomplement of the eigenfunctions (see
\cite{Wy1}); on the other hand any affine extension of that
rotation is spectrally the same, while it cannot have the ELF
property.

Prikhodko and Thouvenot (private communication) have constructed
weakly mixing and non-mixing rank one automorphisms which enjoy
the ELF property.

\section{Inducing and spectral theory}\label{czesc11}
Assume that $T$ is an ergodic automorphism of a standard
probability Borel space $\xbm$. Can ,,all" dynamics be obtained by
inducing (see the glossary) from one fixed automorphism was a
natural question from the very beginning of ergodic theory.
Because of Abramov's formula for entropy $h(T_A)=h(T)/\mu(A)$  it
is clear that positive entropy transformations cannot be obtained
from inducing on a zero entropy automorphism. However here we are
interested in spectral questions and thus we ask how many spectral
types we obtain when we induce. It is proved in \cite{Fr-Or} that
the family of $A\in\cb$ for which $T_A$ is mixing is dense for the
(pseudo) metric $d(A_1,A_2)=\mu(A_1\triangle A_2)$. De la Rue
\cite{Rue1} proves the following result: {\em For each ergodic
transformation $T$ of a standard probability space $\xbm$ the set
of $A\in\cb$ for which the maximal spectral type of $U_{T_A}$ is
Lebesgue is dense in $\cb$.} The multiplicity function is not
determined in that paper.  Recall (without giving a formal
definition, see \cite{Or-Ru-We}) that a zero entropy automorphism
is {\em loosely Bernoulli} (LB for short) if and only if  it can
be induced from an irrational rotation (see also \cite{Feldman},
\cite{Ka1}). The LB theory shows that not all dynamical systems
can be obtained by inducing from an ergodic rotation. However an
open question remained whether LB systems exhaust spectrally all
Koopman representations. An interesting question of M.\ Ratner \cite{Ra0} is whether from every ergodic automorphism $T$ one can induce an automorphism which has countable Lebesgue spectrum (Ratner in  \cite{Ra0} shows that this can be done if $T$ is an irrational rotation).

 In a deep paper \cite{delaRue0}, de la
Rue studies LB property in the class of Gaussian-Kronecker
automorphisms, in particular he constructs $S$ which is not LB.
Suppose now that $T$ is LB and for some $A\in\cb$, $U_{T_A}$ is
isomorphic to $U_S$. Then by the Foia\c{s}-Stratila Theorem, $T_A$
is isomorphic to $S$, and hence $T_A$ is not LB. However, an
induced automorphism from an LB automorphism is LB, a
contradiction.

Another fruitful source of non LB systems comes from taking Cartesian products of some natural LB systems. In \cite{Or-Ru-We}, it is proved  that there exists a rank one (and hence LB) system whose Cartesian square is not LB. Moreover, in \cite{Rat1}, it was shown that the square of the horocycle flow is not LB (the horocycle flow itself being LB, \cite{Ra0}). Recently, \cite{KdlR}, the authors showed that there are staircase rank one transformations whose Cartesian product is not LB.

\section{Rigid sequences} Recall (see the glossary) that an automorphism $T$ of a standard probability Borel space $\xbm$ is called {\em rigid} if there exists a strictly increasing sequence $q_n\to\infty$ such that $\mu(T^{q_n}A\triangle A)\to0$ as $n\to\infty$, for each $A\in\mathcal{B}$.\footnote{In fact, to have a global rigidity sequence, as observed by Thouvenot,  we only need to know that for each $A\in\mathcal{B}$ there is a sequence $(q_{n,A})$ so that $\mu(T^{q_{n,A}}A\triangle A)\to0$.} Equivalently, for each $f\in L^2\xbm$, $U_T^{q_n}f\to f$ in $L^2\xbm$ (it is not hard to see that the latter is equivalent that $\int_0^1e^{2\pi iq_nx}d\sigma_f(x)\to 1$ for any $\sigma_f$ representing the maximal spectral type of $T$, $\|f\|=1$). We call $(q_n)$ a {\em rigidity sequence} of $T$.
Rigidity is one of (purely spectral) fundamental phenomena in ergodic theory. Assuming that $T$ is aperiodic, it is not hard to see that for any rigidity sequence $(q_n)$, we must have $q_{n+1}-q_n\to\infty$. Typical automorphism is rigid and weakly mixing but since weak mixing implies $U^n_T\to0$ weakly on $L^2_0\xbm$ along a sequence of $n$ of full density, there is no much ``room'' left for rigidity sequences. So positive density sequences cannot be rigid but  beyond that, in the class of zero density sequences there can be other, for example algebraic in nature, obstructions for rigidity. For example, as noticed in \cite{Be-Ju-Le-Ro}, \cite{Ei-Gr}, if $P \in\Q[x]$ is any non-zero polynomial taking integer values on $\Z$ then the sequence
$(P(n))$ cannot be rigid for any ergodic automorphism. It is also easy to see that $(2^n)$ is a rigidity sequence while $(2^n+1)$ is not.

A systematic study of sequences which can be rigidity sequences was originated in \cite{Be-Ju-Le-Ro} and \cite{Ei-Gr}. Both papers use harmonic analysis approach to construct rigid sequences (via the standard Gaussian functor). Other constructions are also presented in \cite{Be-Ju-Le-Ro}: rank one constructions, weighted operators, Poisson suspensions, while \cite{Ei-Gr} rather concentrates on so called linear dynamical systems and studies rigidity for weakly mixing automorphisms. One of the results in \cite{Be-Ju-Le-Ro} and \cite{Ei-Gr} states that
if either $q_{n+1}/q_n\to\infty$ or $q_{n+1}/q_n$ is an integer then $(q_n)$ is a rigidity sequence (for a weakly mixing automorphism). On the other hand, Eisner and Grivaux in \cite{Ei-Gr} give an example of a rigid sequence$(q_n)$, for a weakly mixing automorphism,    such that $q_{n+1}/q_n\to1$. As a matter of fact, both \cite{Be-Ju-Le-Ro} and \cite{Ei-Gr} deal with the case of denominators $(q_n)$ of an irrational $\alpha\in[0,1)$ (which are obviously rigidity sequences for the corresponding irrational rotations $Tx=x+\alpha$ on the additive circle) to show that such sequences are rigid for some weakly mixing automorphisms. Let us also mention Aaronson's result \cite{AaRig}: {\em Given any sequence $(r_n)$ of density~0, there is a
sequence $(q_n)$ such that $q_n<r_n$, $n\geq1$, and $(q_n)$  is rigid for some weakly mixing automorphisms.}

Moreover, in \cite{Be-Ju-Le-Ro} two basic questions have been formulated: {\em Given any sequence rigid for some $T$ with discrete spectrum, must it be rigid for some weakly mixing automorphism? What about the converse?}

The positive answer to the first question was given by Adams \cite{Ad1} and Fayad and Thouvenot \cite{Fa-Th}. On the other hand, surprisingly, Fayad and Kanigowski \cite{Fa-Ka} answered negatively the second question: there are rigidity sequences (for weakly mixing automorphisms) which are not rigidity sequences for any rotation. For a strengthening of this result (the existence of a rigidity sequence which, as a subset, is dense in the Bohr topology on $\Z$), see \cite{Gri}.

One can also consider a notion stronger than rigidity, called IP-rigidity (see e.g.\ \cite{Be-Ju-Le-Ro}):
 $(q_n)$ is an IP-rigidity sequence for an automorphism $T$ acting on $\xbm$ if  $T^\xi\to Id$ (in the strong topology of $L^2\xbm$) in the IP sense, that is, when $\xi\to\infty$, where $\xi=q_{m_1}+\ldots+q_{m_k}$ and we require that the smallest element $q_{m_1}$ is going to $\infty$.
 This notion  is studied in \cite{Aa-Ho-Le} relating it to non-singular ergodic theory (more precisely, to groups of so called $L^\infty$-eigenvalues of non-singular automorphisms).
As proved in \cite{Aa-Ho-Le}, in this category, the answer to the second question (above) from \cite{Be-Ju-Le-Ro} turns out to be positive. Moreover, the paper provides an example of a super lacunary sequence (which must be a rigidity sequence by \cite{Be-Ju-Le-Ro} and \cite{Ei-Gr}) which is not an IP-rigid.

In the recent preprint \cite{Ba-Gr-Ma}, rigidity sequences  are compared to other classical notions in harmonic analysis. It is proved that rigidity sequences $(q_n)$ are nullpotent, i.e.\ there exists a topology $\tau$ on $\Z$ making it a topological group such that $q_n\to0$ but they  are never Kazhdan.\footnote{A subset $B\subset\Z$ is called Kazhdan if there exists $\vep>0$ such that each unitary operator $U$ on a separable Hilbert space $H$ having a unit vector $x$ with $\sup_{n\in B}\|U^nx-x\|<\vep$ has a non-zero fixed point.} We find also there a rather surprising result that the family of all rigidity sequences considered as a subset in $\Z^{\N}$ is Borel.

\section{Spectral theory of parabolic dynamical systems}
We say a system is {\em algebraic} if it is a $\Z$ (or $\R$) translation on a quotient of a Lie group by a lattice. Spectral theory (and mixing properties) of {\em algebraic systems} is by now well understood. The two main classes are actions on quotients of semi-simple and nilpotent Lie groups. In the first case, the two main examples are horocycle and geodesic flows on quotients of $SL(2,\R)$.
More generally, one can talk about quasi-unipotent and partially hyperbolic actions. Recall that in the setting of algebraic actions, being quasi-unipotent is equivalent to {\em zero entropy}, while being partially hyperbolic is equivalent to {\em positive entropy}.

It is known that in both cases the spectrum is countable Lebesgue (we refer the reader to \cite{Ka-Th} for a nice description of spectral theory of horocycle and geodesic flows). Actions on quotients of nilpotent Lie groups are also known to have countable Lebesgue spectrum in the orthocomplement of the eigenspace (we refer to \cite{Parry1} for details). Quantitative mixing (and higher order mixing) of algebraic systems is also well understood (we refer the reader to a recent paper, \cite{BEG}, for general results on decay of correlations for algebraic systems on semi-simple Lie groups).

Much less is known in spectral theory of parabolic systems beyond algebraic world. There is no strict definition for a system to be parabolic. However, characteristic features of parabolic systems are: zero-entropy, polynomial orbit growth, strong mixing and equidistribution properties. We describe some classes of (non-algebraic) parabolic systems below.
One of  the main difficulty in studying non-algebraic systems is a lack of many tools from representation theory, which is available in the algebraic setting. Below, we focus on known results and questions in spectral theory of non-algebraic parabolic systems.
\subsection{Time-changes of algebraic systems}
Perhaps the simplest class of non-algebraic parabolic dynamical systems is given by time-changes (or reparametrizations) of algebraic systems. As for algebraic systems, it is natural to consider separately the cases of time-changes of unipotent systems and nilpotent systems. We do it in two paragraphs below.

\paragraph{Time-changes of unipotent systems}
In recent years we witnessed substantial development in understanding the theory of time-changes of unipotent flows. The first (and most studied) case is that of smooth time changes of horocycle flows. Recall first that M.\ Ratner, \cite{RaT}, established measures and joinings rigidity phenomena for time-changes of horocycle flows that are analogous to the Ratner's theory in algebraic setting. In particular in \cite{RaT2} Ratner proved the $H$-property for all (sufficiently smooth time-changes). It is not known if Ratner's joinings and measure rigidity also holds for time-changes of general unipotent flows.

Mixing for smooth time-changes of horocycle flows was established by Marcus \cite{Marc}, generalizing earlier work of Kushnirenko \cite{Ku} who required additionally small derivative in the geodesic direction. A crucial result for the theory is by L.\ Flaminio and G.\ Forni \cite{FlFo}, where the authors classify all invariant distributions and as a consequence show that a typical time-change is not a (measurable) quasi-coboundary (and hence the time-change is not trivially isomorphic to the original flow). A.\ Katok and J.-P. Thouvenot conjectured, \cite{Ka-Th}, that every sufficiently smooth time change of the horocycle flow has countable Lebesgue spectrum. A partial answer to this conjecture was given by G.\ Forni and C.\ Ulcigrai, \cite{FU}, where the authors show that the maximal spectral type of the time-changed flow remains Lebesgue (see also a result of R.\ Tiedra de Aldecoa, \cite{TAl}, where absolute continuity of the spectrum is proven, and \cite{Tiedra1}, \cite{Tiedra2} for further applications of the commutator method in ergodic theory). A full solution of the Katok-Thouvenot conjecture (i.e.\ countable Lebesgue spectrum) was recently given by G.\ Forni, B.\ Fayad and A.\ Kanigowski, \cite{FFK}. Generalizing the approach from \cite{FU}, L.\ Simonelli, \cite{Sim}, showed that the spectrum of smooth time-changes of general unipotent flows remains Lebesgue. It is not known if the multiplicity is infinite, but it seems that the approach from \cite{FFK} has the potential of being applicable in this setting.

Recall that Ratner's work \cite{Ra} allows one to classify joinings between horocycle flows. Recently, there was progress in understanding joinings for time-changes of horocycle flows. In \cite{KLU} (see also a result in \cite{FlFo2}) the authors show that there is a strong dichotomy for two smooth time-changes of horocycle flows: either the time-change functions are cohomologous or the resulting time-changed flows are disjoint.

Even though quantitative mixing for time-changes of unipotent flows is now well understood, not much is known for quantitative higher order correlations:
\begin{question} Is the decay of higher correlations for non-trivial time-changes of horocycle flows (or more generally, unipotent flows) polynomial?
\end{question}
For trivial time-changes, i.e. for the horocycle flow, decay of higher correlations is indeed polynomial by a recent result of M.\ Bj{\"o}rklund, M.\ Einsiedler and A.\ Gorodnik, \cite{BEG} (in fact this applies to general unipotent flows).

\paragraph{Time-changes of nilpotent systems}
Recall that nilpotent flows are never weakly mixing since they always have a non-trivial Kronecker factor. An interesting question is therefore whether one can improve mixing properties of the system by a time-change. The first result in this direction by A.\ Avila, G.\ Forni and C.\ Ulcigrai \cite{AFU} is that there exists a dense set of smooth functions on the {\em Heisenberg nilfmanifold} such that the resulting time-changed Heisenberg flow is {\em mixing}. This result was strenghtened by D.\ Ravotti in \cite{Rav} to quasi-abelian nilflows, and recently to all nilflows by Avila, Forni, Ravotti and Ulcigrai \cite{AFRU}. It is important to mention that the mixing mechanism is non-quantitative. Therefore two questions are natural to ask: what are mixing properties of general time-changes of nilflows and can one obtain some quantitative mixing results. The only case in which
some progress has been recently made is that of time-changes of Heisenberg nilflows. In \cite{FK} the authors show stretched polynomial decay of correlations for smooth time-changes of full measure set of Heisenberg nilflows (parametrized by the frequency of the Kronecker factor). In the case the flow is of {\em bounded type} the authors prove polynomial speed of decay of correlations. Moreover, in \cite{FK2} the authors show that for time-changes of bounded type Heisenberg nilflows, every non-trivial time-change enjoys the $R$-property and as a consequence is {\em mildly mixing}. Moreover, in the above setting, it also follows that every mixing time-change is mixing of all orders.

Mixing and spectral properties of time-changes of general nilflows are poorly understood. In particular, the following questions seems interesting:

\begin{question}  Are all non-trivial smooth time-changes of general nilflows mixing?
\end{question}
The mixing mechanism from \cite{AFU} (see also \cite{Rav} and \cite{AFRU}) is non-quantitative. Therefore the following question seems to be far more challenging:
\begin{question} Does there exist a smooth time change of a general nilflow with AC (Lebesgue) spectrum?
\end{question}

\subsection{Special flows, flows on surfaces, interval exchange
transformations}\label{czesc12}
In this section we will describe spectral results for special flows over interval exchange transformations (IETs) (irrational rotations begin a particular case). As described below, such special flows arise as representation of smooth  locally hamiltonian flows on surfaces.

\subsubsection{Interval exchange transformations}
To define an interval exchange transformation (IET) of $m$ intervals we need a
permutation $\pi$ of $\{1,\ldots,m\}$ and a probability vector
$\la=(\la_1,\ldots,\la_m)$ (with positive entries). Then we define
$T=T_{\la,\pi}$ of $[0,1)$ by putting
$$
T_{\la,\pi}(x)=x+\beta_i^{\pi}-\beta_i\;\;\mbox{for}\;\;
x\in[\beta_i,\beta_{i+1}),$$ where $\beta_i=\sum_{j<i}\la_j$,
$\beta^{\pi}_i=\sum_{\pi j<\pi i}\beta_j$. Obviously, each IET
preserves Lebesgue measure.  One of possible approaches to study
ergodic properties of IET is an ``a.e" approach ``seen" in the
definition of $T_{\la,\pi}$. It is based on the fundamental fact
that the induced transformation on a subinterval of $[0,1)$ is
also IET (see \cite{Co-Fo-Si}). This leads to a very delicate and
deep mathematics based on Rauzy induction, which is a way of
inducing on special intervals, considering only irreducible
permutations whose set is partitioned into orbits of some  maps
(any such an orbit is called a {\em Rauzy class}). If now $\crr$
is a Rauzy class of permutations and $\la$ lies in the standard
simplex $\Delta_{m-1}$ then the Rauzy induction together with a
natural renormalization leads to a map
$\cp:\crr\times\Delta_{m-1}\to \crr\times\Delta_{m-1}$. For a
better understanding of the dynamics of the Rauzy map Veech
\cite{Ve1} introduced the space of {\em zippered rectangles}. A
zippered rectangle associated to the Rauzy class $\mathcal{R}$ is
a quadruple $(\lambda, h, a, \pi)$, where $\lambda\in
\mathbb{R}^m_+$, $h\in \mathbb{R}^m_+$, $a\in \mathbb{R}^m_+$,
$\pi\in\mathcal{R}$ and the vectors $h$ and $a$ satisfy some
equations and inequalities. Every zippered rectangle $(\lambda, h,
a, \pi)$ determines a Riemann structure on a compact connected
surface. Denote by $\Omega(\mathcal{R})$ the space of all zippered
rectangles, corresponding to a given Rauzy class $\mathcal{R}$ and
satisfying the condition $\langle\lambda,h\rangle=1$. In
\cite{Ve1}, Veech defined a flow $(P^t)_{t\in\R}$ on the space
$\Omega(\mathcal{R})$ putting
\[ P^t(\lambda,h,a,\pi) = (e^t\lambda, e^{-t}h, e^{-t}a, \pi)\]
and extended the Rauzy map. On so called  {\em Veech moduli space}
of zippered rectangles, the flow $(P^t)$ becomes the {\em
Teichm\"uller flow} and it preserves a natural Lebesgue measure
class; by passing to a transversal  and projecting the measure on
the space of IETs $\mathcal{R}\times \Delta_{m-1}$ Veech has
proved the following fundamental theorem (\cite{Ve1}, see also
\cite{Ma}) which is a generalization of the fact that Gauss
measure $\frac1{\ln2}\frac1{1+x}dx$ is invariant for the Gauss map
which sends $t\in(0,1)$ into the fractional part of its inverse.

\begin{Th} [Veech, Masur, 1982] Assume that $\crr$ is a Rauzy
class. There exists a $\sigma$-finite measure $\mu_{\crr}$ on
$\crr\times\Delta_{m-1}$ which is $\cp$-invariant, equivalent to
``Lebesgue" measure, conservative and ergodic.
\end{Th}

In \cite{Ve1} it is proved that a.e.\ (in the above sense) IET is
then of rank one (and hence is ergodic and has a simple spectrum).
He also formulated as an open problem whether we can achieve the
weak mixing property a.e. This has been recently answered in
positive by A. Avila and  G. Forni \cite{Av-Fo} (for $\pi$ which
is not a rotation).

Katok \cite{Ka7} proved that IET  have no mixing factors (in fact
his proof shows more:  the IET class is disjoint from the class of
mixing transformations).  By their nature, IET transformations are
of finite rank (see \cite{Co-Fo-Si}) so they are of finite
multiplicity. They need not be of simple spectrum (see remarks in
\cite{Ka-Th} pp. 88-90). It remains an open question whether  an
IET can have a non-singular spectrum. Joining properties in the
class of exchange of~3 and more intervals are studied in
\cite{Fe-1}, \cite{Fe-2}. In \cite{CHE} the authors show that a.e. $3$-IET is not simple. This answers a special case of a question of Veech
\cite{Ve0} whether a.e.\ IET is simple. The case of $d$-IET's with $d>3$ is still widely open.

\subsubsection{Smooth flows on surfaces and their special representations}
We consider a closed, connected, smooth and orientable surface $S$ of genus $g\geq 1$. Let $X: S\mapsto TS$ be a smooth vector field with finitely many fixed points and such that the corresponding flow $(\phi_t^X)$ preserves a smooth area form $\omega$. The flow $(\phi_t^X)$ is called a {\em locally Hamiltonian flow}; it is locally given by a  smooth hamiltonian $H$ (up to an additive constant), so that $(\phi_t^X)$ is a solution to
\beq\label{flow2}
\frac{dx}{dt}=\frac{\partial H}{\partial y},\;
\frac{dy}{dt}=\frac{-\partial H}{\partial x}.\eeq

There has recently been some progress in understanding ergodic and spectral properties of locally Hamiltonian flows. It follows by \cite{Ar} that $(\phi_t^X)$ can be decomposed into finite number of topological discs filled with periodic orbits and finite number of connected components on which the flow is minimal. We are always interested in properties of $(\phi_t^X)$ on the minimal components. A classical way of studying
locally hamiltonian flows is through their {\em special representations}. On each minimal component one chooses a smooth transversal and represents $(\phi_t^X)$ by the first return map (the base transformation $T$) and the first return time (the roof function $f$). We will denote the special representation of $(\phi_t^X)$ by $T^f$. It is well known that $T$ is either a rotation (in genus $1$ case) or an interval exchange transformation (in higher genus) and the roof function is smooth except finitely many points at which it explodes according to the nature of the fixed points. Notice first that if $(\phi_t^X)$ has no fixed points then $S$ is a two-dimensional torus (by Gauss-Bonet). In this case $(\phi_t^X)$ is a smooth time-change of a linear flow on the two-torus. Therefore it is never mixing \cite{Ko1}, disjoint from mixing flows \cite{Fr-Le} and in particular has purely singular spectrum. From now on we will therefore assume that the set of fixed points is non-empty. If all fixed points are non-degenerated (i.e. hessian of $H$ is non-zero at every critical point) the roof function $f$ has {\em logarithmic singularities}, i.e. it blows up logarithmically at finitely many points. Logarithmic singularities can be either {\em symmetric} (for instance, if there are no saddle connections) or {\em asymmetric} (in case there is a saddle loop). Ergodic properties of $(\phi_t^X)$ depend strongly on whether the roof has symmetric or asymmetric singularities. If there are {\em degenerated} fixed points, Kochergin, \cite{Ko3}, constructed smooth hamiltonian for which the roof function at singularities blows up power-like (like $x^{\gamma}$, $\gamma\in (-1,0)$). By changing a speed as it is
done in \cite{Fr-Le2}  so that critical points of the vector field
in~(\ref{flow2}) become singular points, Arnold's special
representation is transformed to a special flow over the same
irrational rotation however under a piecewise smooth function. If
the sum of jumps is not zero then in fact we come back to von
Neumann's class of special flows considered in \cite{vonNeumann}. Similar classes of special flows (when the roof function is of
bounded variation) are obtained from ergodic components of flows
associated to billiards in convex polygons with rational angles
\cite{Ka-Ze}.

It is convenient to express ergodic and spectral properties of locally hamiltonian flows in terms of their special representation. We will do it in the next subsection.

\subsection{Special flows over rotations and interval exchange transformations}
 Special flows were introduced to ergodic theory by von
Neumann in his fundamental work \cite{vonNeumann} in 1932. Also in
that work he explains how to compute eigenvalues for special
flows, namely: {\em $r\in\R$ is an eigenvalue of $T^f$ if and only
if the following functional equation
$$ e^{2\pi irf(x)}=\xi(x)/\xi(Tx)$$
has a measurable solution $\xi:X\to\T$}. We recall also that the
classical Ambrose-Kakutani theorem asserts that practically each
ergodic flow has a special representation (\cite{Co-Fo-Si},
 see also  Rudolph's theorem on special representation therein). As described in the previous subsection, we will consider the case when the roof functions is: of bounded variation, has logarithmic singularities (symmetric or asymmetric) or has power singularities.

\paragraph{Bounded variation roof function.}
 Kochergin \cite{Ko1} showed that special flows over
irrational rotations and under bounded variation functions are
never mixing. This has been recently strengthened in \cite{Fr-Le}
to the following: {\em If $T$ is an irrational rotation and $f$ is
of bounded variation then the special flow $T^f$ is spectrally
disjoint from all mixing flows}. In particular all such flows have
singular spectra. Moreover, in \cite{Fr-Le} it is proved that
whenever the Fourier transform of the roof function $f$ is of
order O$(\frac1n)$ then $T^f$ is disjoint from all mixing flows
(see also \cite{Fr-Le3}). In fact in the papers \cite{Fr-Le} --
\cite{Fr-Le2} the authors discuss the problem of disjointness of
those special flows with all ELF-flows conjecturing that no flow
of probabilistic origin has a smooth realization on a surface. In
\cite{Le-Wy}, the analytic case is considered leading to a ``generic"
result on disjointness with the ELF class generalizing the
classical Shklover's result on the weak mixing property \cite{Sh}. A. Katok  \cite{Ka7} proved the absence of mixing for special
flows over IET when the roof function is of bounded variation (see
also \cite{Ry12}). Katok's theorem was strengthened in
\cite{Fr-Le4} to the disjointness theorem with the class of mixing
flows.  A.\ Avila and G.\ Forni  \cite{Av-Fo} proved that a.e.\ translation
flow on a surface (of genus at least two) is weakly mixing (which
is a drastic difference with the linear flow case of the torus,
where the spectrum is always discrete).

One important property in the class of special flows over rotations and IET's is Ratner's property (R-property). This property may be viewed as a particular
way of divergence of orbits of close points; it was shown to hold
for horocycle flows by M. Ratner  \cite{Ra}. We refer the reader
to \cite{Ra} and the survey article \cite{Th} for the formal
definitions and basic consequences of R-property. In particular,
R-property implies ``rigidity" of joinings and it also implies the PID
property; hence mixing and R-property imply mixing of all orders.
In \cite{Fr-Le2}, \cite{Fr-Le-Le} a version of R-property is shown
for the class of von Neumann special flows (however $\alpha$ is
assumed to have bounded partial quotients). This allowed one to
prove there that such flows are even mildly mixing (mixing is
excluded by a Kochergin's result). The eigenvalue problem (mainly how many frequencies can have the
group of eigenvalues) for special flows over irrational rotations
is studied in \cite{Fa-Ka-Wi}, \cite{Fa-Wi}, \cite{Gu-Pa}.

It follows from \cite{Fr-Le} that von Neumann flows have singular spectrum. However nothing is known about their multiplicity.
\begin{question}
What is the spectral multiplicity of von Neumann flows?
\end{question}

\paragraph{Symmetric logarithmic singularities.}
Kochergin \cite{Ko2} proved the absence of mixing  for flows where
the roof function has finitely many singularities, however some Diophantine restriction is put on $\alpha$. In \cite{Le10}, where also the absence of mixing is considered for the symmetric logarithmic case,  it was conjectured (and proved
for arbitrary rotation) that a necessary condition for mixing of a
special flow $T^f$ (with arbitrary $T$ and $f$) is the condition
that the sequence of distributions $((f_0^{(n)})_\ast)_n$ tends to
$\delta_\infty$ in the space of probability measures on $\ov{\R}$.
K. Schmidt \cite{Sch0} proved it using the theory of cocycles and
extending a result from \cite{Aa-We} on tightness of cocycles. Ulcigrai, \cite{uczennicaS}, showed that for every $d\geq 2$ and a.e.\ IET of $d$ intervals, the corresponding special flow $T^f$ is {\em not} mixing. However, Chaika and Wright, \cite{CWr}, proved existence of an IET $T$ such that $T^f$ is mixing. Notice that by \cite{Fr-Le} it also follows that for a.e.\ irrational rotation $T^f$ has purely singular spectral type. Recent result, \cite{CFKU}, shows that one can also prove singularity of the spectrum for {\em symmetric IET's} in the base. This in particular shows that minimal flows on genus $2$ surfaces (with two isometric saddles) have purely singular spectral type. In \cite{KaKu} the authors showed that if $T$ is an IET of {\em bounded type} then $T^f$ is mildly mixing (and has the $R$-property). For symmetric logarithmic singularities the following two questions are open:
\begin{question}What is the maximal spectral type of $T^f$ for a general IET's?
\end{question}

Nothing is known about multiplicity of the spectrum in this setting.

\paragraph{Asymmetric logarithmic singularities.}
In this case mixing properties are different. Khanin and Sinai, \cite{Kh-Si}, showed that $T^f$ is mixing for a.e.\ irrational rotation $T$. This was strengthened by Ulcigrai to a.e. IET, \cite{uczennicaS2}. Moreover, \cite{Rav22}, Ravotti obtained quantiative mixing estimates (with sub-logarithmic speed of decay of correlations). Not much was known about multiple mixing for $T^f$. This changed recently: in \cite{Fak} the authors showed that for a.e. irrational rotation,  $T^f$ enjoys a variant of the $R$-property and hence is multiple mixing. This was strenghtened to a.e. IET's in \cite{KKU}. The following two questions seem to be natural (the first one already stated as Questions 34, 35 in \cite{FKri}):

\begin{question}What is the maximal spectral type of $T^f$?
\end{question}
Moreover, one can ask about quantitative higher order decay:
\begin{question}Is the decay of higher order correlations sub-logarithmic?
\end{question}
Both of the above question are open even for rotations in the base.
As in the symmetric case, nothing is known about multiplicity of the spectrum.

\paragraph{Power singularities.}
In case $f$ has power singularities, $T^f$ was shown to be mixing by Kochergin, \cite{Ko3} (for any uniquely ergodic irrational rotation and IET's). In \cite{Fak} the authors show that if $T$ is a rotation of bounded type, then $T^f$ is multiple mixing (it enjoys a variant of the $R$-property). Polynomial decay for some flows with power singularities was obtained by B.\ Fayad, \cite{Fa4}. In \cite{FFK} the authors considered the spectrum of $T^f$. They showed that if $f$ has sufficiently strong power singularity (of the form $x^{-1+\eta}$ for small $\eta>0$) then $T^f$ has countable Lebesgue spectrum for a.e.\ irrational rotation. To the best of our knowledge, this is the only result dealing with multiplicity for smooth surface flows. The following problems are natural (see Question 34 in \cite{FKri}):
\begin{question}
What is the maximal spectral type of $T^f$ when $f$ has power singularities? What is the multiplicity?
\end{question}
To answer this one needs to consider general IET's in the base, as well as functions with weaker power singularity than in \cite{FFK}. The following question is still open (Question 38 in \cite{FKri}):
\begin{question}Are all mixing surface flows mixing of all orders?
\end{question}

Finally, it  may also be useful to show that smooth flows
on surfaces are disjoint from flows of probabilistic origin -- see
\cite{Ju-Le2}, \cite{Ju-Le3}, \cite{Le-Pa-Ro}, \cite{Ry-Th},
\cite{Th1}.

B. Fayad \cite{Fa2} gives a criterion that implies singularity of
the maximal spectral type for a dynamical system on a Riemannian
manifold. As an application he gives a class of smooth mixing
flows (with singular spectra) on $\T^3$ obtained from linear flows
by a time change (again this is a drastic difference with
dimension two, where a smooth time change of a linear flow leads
to non-mixing flows \cite{Co-Fo-Si}).

We mention at the end that if we drop here (and in other problems)
the assumption of regularity of $f$ then the answers will be
always positive because of the LB theory; in particular, there is a
section of any horocycle flow (it has the LB property \cite{Ra0})
such that in the corresponding special representation  $T^f$, the map
$T$ is an irrational rotation. Using a Kochergin's result
\cite{Ko9} on cohomology (see  also \cite{Ka0}, \cite{Ru4}) the
$L^1$-function $f$ is cohomologous to a positive function $g$
which is even continuous, thus $T^f$ is isomorphic to $T^g$.

\section{ Spectral theory for locally compact groups of type $I$}
\label{nieab}
This section has been written by A.\ Danilenko.
\subsection{Groups of type $I$}
The spectral theory presented here for Abelian group actions extends potentially  to probability preserving actions of non-Abelian locally compact groups of type $I$.
We now provide  the  definition of type $I$.
Let $G$ be a locally compact second countable group, $\mathcal{H}$ a separable  Hilbert space
and $\pi:G\ni g\mapsto \pi(g)$ is a (weakly) continous unitary representation of $G$ in $\mathcal{H}$.
We say that $\pi$ is {\it of type $I$} if there is a subset $A\subset\{1,2,\dots,+\infty\}$ such that
$\pi$ it is unitarily equivalent to the orthogonal sum
$
\bigoplus_{k\in A} U_k\otimes I_k,
$
where $U_k$ is  a unitary representation of $G$ with a simple spectrum and $I_k$ is the trivial representation in the Hilbert space of dimension $k$.

\begin{Def}\em
If every unitary representation of $G$ is of type $I$ then $G$ is called {\it of type $I$}.
\end{Def}

Denote by $\widehat G$ the {\it unitary dual} of $G$, i.e.\ the set of  unitarily equivalent classes of all irreducible unitary representations of $G$.
If $G$ is Abelian then every irreducible representation of $G$ is 1-dimensional.
Hence $\widehat G$ is identified naturally with the group of characters of $G$.
In the general case, let Irr$_n(G)$ stand for the set of all irreducible unitary representations of $G$ in the $n$-dimensional separable Hilbert space $\mathcal{K}_n$, $1\le n\le+\infty$.
Endow it with the natural Borel structure, i.e.\ the smallest one in which the mapping
$\pi\mapsto\langle\pi(g)f,h\rangle$ is Borel for every $g\in G$ and $f,h\in\mathcal{K}_n$.
It is standard.
Let $\widehat G_n$ is the quotient of Irr$_n(G)$ by the unitary equivalence relation.
Endow $\widehat G_n$ with the quotient Borel structure.
Since $\widehat G=\bigsqcup_{n=1}^\infty \widehat G_n\sqcup \widehat G_\infty$,
we obtain a Borel structure on $\widehat G$.
It is called {\it Mackey Borel structure} on $\widehat G$.
By the Glimm theorem, the Mackey Borel structure is standard if and only if  $G$ is of type $I$
\cite{Gli}.

For  $n\in\Bbb N\cup\{\infty\}$, denote by $I_n$ the identity operator on $\mathcal{K}_n$.
Then  for each unitary representation $\pi$ of $G$ in $\mathcal{H}$, there is a measure $\lambda$,
a measurable field $\widetilde G\ni\omega\mapsto\mathcal{H}_\omega$ of Hilbert spaces,
a measurable field
 $\widehat G\ni\omega\mapsto V_\omega$ of irreducible unitary $G$-representations
such that $V_\omega\in\omega$ and $\mathcal{H}_\omega$ is the space of $V_\omega$
on $\widehat G$ and a measurable map $m:\widehat G\to\Bbb N\cup\{+\infty\}$ such that
$$
\mathcal{H}=\int^{\oplus}_{\widehat G}\mathcal{H}_w\otimes \mathcal{K}_{m(w)}\,d\lambda(\omega)\;\mbox{and}$$$$
\pi(g)=\int^{\oplus}_{\widehat G} V_w(g)\otimes I_{m(w)}\,d\lambda(\omega).
$$
It appears that if $G$ is of type $I$ then the equivalence class of $\lambda$ is defined uniquely by $\pi$ and the function $m$ is defined up to a $\lambda$-zero subset.
We  call the class of $\lambda$ {\it the maximal spectral type of $\pi$}, and we call $m$ {\it the spectral multiplicity of $\pi$}.
If we have a probability preserving action $T=(T_g)_{g\in G}$ of $G$ then we can consider the corresponding Koopman unitary representation $\pi$ of $G$.
The maximal spectral type of $\pi$ and the spectral multiplicity of $\pi$ is called {\it the maximal spectral type  of $T$} and {\it the spectral multiplicity of $T$} respectively.
If $G$ is Abelian then these concepts coincide with their classic counterparts considered above in the survey.

All compact groups, Abelian groups, connected semisimple Lie groups, nilpotent Lie groups (or, more generally, exponential Lie groups) are of type $I$.
Each subgroup of $GL(n,\Bbb R)$ determined by a system of algebraic equations is also of type $I$.
Solvable Lie groups can be as of type $I$ as not of type $I$.
If $G$ is a   countable (discrete) groups then $G$ is of type $I$ if and only if it is virtually Abelian, i.e.\ it contains an Abelian subgroup of finite index \cite{Thom}.

Even if we know that a non-Abelian group $G$ is of type $I$, it is usually not an easy problem to describe  $\widehat G$ explicitly.
Kirillov introduced an {\it orbit method} for a description of $\widehat G$ when $G$ is a
connected, simply connected nilpotent Lie group (the method was developed further for solvable groups).
He identified $\widehat G$ with the space of  orbits for the co-adjoint $G$-action on the dual space $\mathfrak{g}^*$ of its Lie algebra $\mathfrak{g}$.
Though Kirillov's method gives an algorithm how to describe $\widehat G$,  not so many groups are known  for which the unitary dual is described explicitly.

\subsection{Spectral properties of Heisenberg group actions}
The  3-dimensional Heisenberg group $H_3(\Bbb R)$ is perhaps
one of the simplest  examples of non-Abelian nilpotent Lie groups for which the orbit method leads to a very  concrete description of  the unitary dual \cite{Kir}.
Recall that
$$
H_3(\Bbb R)=\left\{\left(\begin{array}{ccc}
1& a & c\\
0& 1 & b\\
0& 0 & 1
\end{array}\right)
: a,b,c\in\Bbb R
\right\}.
$$
For simplicity, we will denote the matrix $\left(\begin{array}{ccc}
1& a & c\\
0& 1 & b\\
0& 0 & 1
\end{array}\right)$ by $[a,b,c]$.
The unitary dual  $\widehat{H_3(\Bbb R)}$ is identified with $\Bbb R^2\sqcup\Bbb R^*$ endowed with the natural standard Borel structure.
Every irreducible unitary representation of $H_3(\Bbb R)$ is unitarily equivalent to either a one-dimensional $\pi_{\alpha,\beta}$ with $(\alpha,\beta)\in\Bbb R^2$ or an  infinite dimensional
$\pi_\gamma$ in $L^2(\Bbb R,\mbox{Leb})$, with $\gamma\in\Bbb R^*:=\Bbb R\setminus\{0\}$, such that
$$
\pi_{\alpha,\beta}[a,b,c]:=
e^{2\pi i(\alpha a+\beta b)},
$$$$
\pi_{\gamma}[a,b,c]
f(x):=e^{2\pi i(c+bx)}f(x+a),\quad f\in L^2(\Bbb R,\mbox{Leb}).
$$
Now,  given a measure preserving $H_3(\Bbb R)$-action on a standard probability space
$(X,\mu)$, let $U$ denote the corresponding Koopman representation in $L^2(X,\mu)$.
Then there are a  probability measure $\sigma^{1,2}$ on $\Bbb R^2$, a function
$l^{1,2}:\Bbb R^2\to\Bbb N$ a probability measure $\sigma^3$ on $\Bbb R^*$, a function
$l^3:\Bbb R^*\to\Bbb N$ such that
$$
L^2(X,\mu)=
\int_{\Bbb R^2}^\oplus\,\bigoplus_{j=1}^{l^{1,2}(\alpha,\beta)}\Bbb C\,d\sigma^{1,2}(\alpha,\beta)\oplus
\int_{\Bbb R^*}^\oplus\,\bigoplus_{j=1}^{l^{3}(\alpha,\beta)}L^2(\Bbb R,\mbox{Leb})\,d\sigma^{3}(\gamma)\; \mbox{and}
$$$$
U=\int_{\Bbb R^2}^\oplus\,\bigoplus_{j=1}^{l^{1,2}(\alpha,\beta)}\pi_{\alpha,\beta}\,d\sigma^{1,2}(\alpha,\beta)\oplus
\int_{\Bbb R^*}^\oplus\,\bigoplus_{j=1}^{l^{3}
(\alpha,\beta)}\pi_{\gamma}\,d\sigma^{3}(\gamma).
$$
We  now compare the spectral properties of $T$ with the spectral properties
of the restriction $T$ to the center of $H_3(\Bbb R)$.
The center is the subgroup $\{[0,0,t]:t\in\Bbb R\}$.

\begin{Prop} [\cite{Dan}]
The maximal spectral type of the flow $(T_{[0,0,t])})_{t\in\Bbb R}$
contains the measure
$\sigma^{1,2}(\Bbb R^2)\delta_0+\sigma^3$.
 The corresponding spectral multiplicity is
 $ \int_{\Bbb R^2}l^{1,2}d\sigma^{1,2}$  at the point $t=0$
 and  the infinity if $t\in\Bbb R^*$.
\end{Prop}

\begin{Th} [\cite{Dan}]
If  $(T_{[0,0,t])})_{t\in\Bbb R}$ is ergodic  then:
\begin{enumerate}
 \item
$\sigma^{1,2}(\Bbb R^2\setminus\{(0,0)\}) = 0$, i.e.\ there are no non-trivial one-dimensional representations in the spectral decomposition of $U$.
The maximal spectral type of $T$ equals the maximal spectral type of the restriction of $T$
to the center of $H_3(\Bbb R)$ (modulo the natural identification);
\item
 $T$ is mixing (see also \cite{Ry1}, \cite{Ry2});
\item
the weak closure of the group $\{T_g:g\in H_3(\Bbb R)\}$ in Aut$(X, \mu)$
is the union of $\{T_g:g\in H_3(\Bbb R)\}$ and the weak closure of
$\{T_{[0,0,t]}:t\in\Bbb R\}$;
\item
if $T$ is rigid then $(T_{[0,0,t])})_{t\in\Bbb R}$ is rigid.
\end{enumerate}
\end{Th}

In \cite{Dan}, there were constructed explicit  examples of  mixing of all orders rank-one (and hence zero entropy) actions
$T$ of the Heisenberg group.

The concept of simplicity  for ergodic $H_3(\Bbb R)$-actions is defined
 in a similar way as for the Abelian  actions.
 A simple $H_3(\Bbb R)$-action  $T$ {\it has MSJ} if the centralizer of the action is $\{T_{[0,0,t]}: t\in\Bbb R\}$.
It was  shown  in \cite{Dan} that the  examples of mixing $H_3(\Bbb R)$-actions constructed there satisfy also the following:
\begin{enumerate}
\item
The flow $(T_{[0,0,t])})_{t\in\Bbb R}$ is simple and the centralizer of it equals the group
$\{T_g:g\in H_3(\Bbb R)\}$.
\item
The transformation $T_{[0,0,1]}$ is simple and the centralizer of it equals
$\{T_g:g\in H_3(\Bbb R)\}$.
\item
$T$ has MSJ.
\end{enumerate}
As a corollary, we obtain examples of  mixing Poisson and mixing Gaussian (probability preserving) actions of $H_3(\Bbb R)$.

\subsection{Heisenberg odometers}
In \cite{DanLem}, the authors isolated a special  class of ergodic $H_3(\Bbb R)$-actions, called {\it the odometer} actions.
Namely, let $\Gamma_1\supset\Gamma_2\supset\cdots$ be a sequence of lattices in $H_3(\Bbb R)$.
Then, we can associate a sequence of homogeneous $H_3(\Bbb R)$-spaces intertwined with
$H_3(\Bbb R)$-equivariant maps:
$$
H_3(\Bbb R)/\Gamma_1\longleftarrow H_3(\Bbb R)/\Gamma_2\longleftarrow\cdots.
$$
Denote by $X$ the projective limit of this sequence.
Then $X$ is a compact Polish $G$-space.
Endow each space  $H_3(\Bbb R)/\Gamma_n$ with the  Haar probability measure.
The projective limit of the sequence of these measures is a $H_3(\Bbb R)$-invariant probability measure $\mu$.
Of course, the $H_3(\Bbb R)$-action on $(X,\mu)$ is ergodic.
It is called {\it the Heisenberg odometer  associated with}  $(\Gamma_n)_{n=1}^\infty$.
We consider the Heisenberg odometers as  non-commutative counterparts of the
ergodic $\Bbb Z$-actions and $\Bbb R$-actions with pure point rational spectrum.

A complete spectral decomposition of the Heisenberg odometers is found in \cite{DanLem}.
Denote by $p:H_3(\Bbb R)\to\Bbb R^2$ the homomorphism
$
[a,b,c]\mapsto(a,b).
$
The kernel of this homomorphism is the center of the Heisenberg group.
Given a lattice $\Gamma$ in $H_3(\Bbb R)$, we denote by $\xi_\Gamma$ a positive real
such that $\Gamma\cap\mbox{Ker}\,p=\{[0,0,n\xi_\Gamma]\mid n\in\Bbb Z\}$.
Of course, $p(\Gamma)$ is a lattice in $\Bbb R^2$.
If $p(\Gamma) =A(\Bbb Z^2)$ for some matrix $A\in GL(2,\Bbb R)$ then we denote by $p(\Gamma)^*$ the dual lattice $(A^*)^{-1}\Bbb Z^2$ in $\Bbb R^2$.

\begin{Th} Let $U$ stand for the Koopman unitary representation of the Heisenberg odometer associated with a sequence of lattices $\Gamma_1\supset\Gamma_2\supset\cdots$.
If $\bigcup_{n=1}^\infty p(\Gamma_n)^*$ is not closed in $\Bbb R^2$ then
$$
U=\bigoplus_{(\alpha,\beta)\in\bigcup_{n=1}^\infty p(\Gamma_n)^*}\pi_{\alpha,\beta}\oplus
\bigoplus_{0\ne\gamma\in\bigcup_{n=1}^\infty\xi_{\Gamma_n}^{-1}\Bbb Z} \bigoplus_1^\infty\pi_\gamma.
$$
\end{Th}

An analogous decomposition is found also for the case where $\bigcup_{n=1}^\infty p(\Gamma_n)^*$ is  closed in $\Bbb Z^2$ (see \cite{DanLem} for details).
Thus, we see that the maximal spectral type of Heisenberg  odometers  is  purely atomic.

For a decreasing sequence $\Gamma=(\Gamma_n)_{n=1}^\infty$
of lattices in $H_3(\Bbb R)$ we let $S(\Gamma):=\bigcup_{n=1}^\infty p(\Gamma_n)^*$
and $\xi_\Gamma:=\bigcup_{n=1}^\infty\xi_{\Gamma_n}^{-1}\Bbb Z$.
The following theorem (except for the first claim) and the below remarks demonstrate a drastic difference between $H_3(\Bbb R)$-odometers
and $\Bbb Z$-odometers.

\begin{Th} Two Heisenberg odometers $T$ and $T$ associated with decreasing sequences of latices $\Gamma$ and $\Gamma'$ respectively are unitarily equivalent if and only if $S(\Gamma)=S(\Gamma')$ and $\xi_\Gamma=\xi_\Gamma'$.
The direct product $T\times T':=(T_g\times T_g')_{g\in G}$ is not spectrally equivalent to any Heisenberg odometer.
$T\times T' $ is ergodic if and only if $S(\Gamma)\cap S(\Gamma')=\{0\}$.
$T\times T' $ is ergodic and has discrete maximal spectral type if and only if $S(\Gamma)\cap S(\Gamma')=\{0\}$ and $\xi_\Gamma\cap \xi_\Gamma'=\{0\}$.
\end{Th}

It was shown in \cite{DanLem} that   Heisenberg odometers are not {\it isospectral}, i.e.
the unitary equivalence, in general, does not imply isomorphism for the underlying $H_3(\Bbb R)$-actions.
It was also shown in \cite{DanLem} that  Heisenberg odometers are not {\it spectrally determined}:
 a Heisenberg odometer is constructed which is unitarily equivalent to an $H_3(\Bbb R)$-action which is not isomorphic to any Heisenberg odometer.

\subsection{On  the ``finitely dimensional'' part of the spectrum}
 Suppose  that $G$ is an arbitrary locally compact second countable group.
 If $G$ is not of type $I$ then  $\widehat G$  furnished with the Mackey Borel structure is a ``bad'' (not standard) Borel space.
 Then a decomposition of a unitary representation of $G$ into irreducibles can be done in essentially non-unique way (see \cite{Kir}).
 Nevertheless, this ``badness'' is related only to the infinite-dimensional part of the spectrum.
 Thus, the union $\bigsqcup_{n=1}^\infty \widehat G_n=\widehat G\setminus\widehat G_\infty$ is a "good" standard Borel space.
Thus, given a measure preserving $G$-action, we can  study ``the finitely dimensional part" of the spectrum of the corresponding Koopman unitary representation of $G$.

This approach was used by Mackey in \cite{Mac}, where he made an attempt to extend the  theory of action with pure discrete spectrum to non-Abelian groups.
For that he isolated a class of ergodic $G$-actions $T$ for which the Koopman representation $U_T$ decomposes into a (countable)  family of  finite dimensional irreducible  representations.
We note that the family is uniquely defined by $T$.
Mackey called such $T$ {\it an action with pure point spectrum}.
He established a structure for these actions.
He showed that  $T$ has a pure point spectrum if and only if $T$ is isomorphic to a $G$-action by rotations on a homogeneous space of $G$ by a compact subgroup.
  However, in general, in contrast with the Abelian case, the $G$-actions with pure point spectrum are not necessarily  isospectral even in the case of finite $G$
(see \cite{Tod} and  \cite{Le--We} Section~6 for  counterexamples).

In \cite{Li--Ug}, Lightwood, {\c S}ahin, and  Ugarcovici considered certain
odometer actions of the discrete Heisenberg group $H_3(\Bbb Z):=\{[a,b,c]: a,b,c\in\Bbb Z\}$.
This group is not of type $I$ (in contrast with $H_3(\Bbb R)$) because each subgroup of finite index in $H_3(\Bbb Z)$  is non-Abelian.
Hence $\widehat{H_3(\Bbb Z)}$ is a ``bad'' Borel space.
On the other hand,  the  (standard Borel) subspace $\bigsqcup_{n=1}^\infty\widehat{H_3(\Bbb Z)}_n$ of it is explicitly described in \cite{In}.
Consider now an  $H_3(\Bbb Z)$-odometer $T$ generated by a decreasing sequence
$\Gamma_1\supset\Gamma_2\supset\cdots$ of normal subgroups of finite index in $H_3(\Bbb Z)$.
We call such an odometer {\it normal}.
 Then  $T$ has a pure point  spectrum in the sense of Mackey.
Moreover,  the full list  of the finitely dimensional irreducible unitary representations that occur in  $U_T$
is found in \cite{Li--Ug} in terms of the sequence $(\Gamma_n)_{n=1}^\infty$.
It was shown later in \cite{DanLem} that the multiplicity of each irreducible component in $U_T$ equals the dimension of this component.
It was also proved  in \cite{DanLem} that the normal $H_3(\Bbb Z)$-odometers are isospectral.
Thus, the unitary equivalence of the Koopman representations implies isomorphism of the underlying normal $H_3(\Bbb Z)$-odometers.

\section{Future directions}\label{czesc13}
We have already seen several cases where spectral properties
interact with measure-theoretic properties of a system.
Let us mention a few more cases which require further research and
deeper understanding.

We recall that the weak mixing property can be understood as a
property complementary to discrete spectrum (more precisely to the
distality \cite{Fu2}), or similarly mild mixing property is
complementary  to rigidity. This can be phrased quite precisely by
saying that $T$ is not weakly (mildly) mixing if and only if it
has  a non-trivial factor with discrete spectrum (it has a
non-trivial rigid factor). It has been a question for quite a long
time if in a sense mixing can be ``built" on the same principle.
In other words we seek a certain ``highly" non-mixing factor. It
was quite surprising when in 2005 F. Parrreau (private
communication) gave the positive answer to this problem.

\begin{Th}[Parreau] Assume that $T$ is an ergodic automorphism of
a standard probability space $\xbm$. Assume moreover that $T$ is
not mixing. Then there exists a non-trivial factor (see below) of
$T$ which is disjoint  from all mixing automorphisms.
\end{Th}

In fact, Parreau proved that each factor of $T$ given by
$\cb_\infty(\rho)$ (this $\sigma$-algebra is described in
\cite{Le-Pa-Th}), where $U_{T}^{n_k}\to J_\rho$, is disjoint from
all mixing transformations. This proof leads to some other results
of the same type, for example: {\em Assume that $T$ is an ergodic
automorphism of a standard probability space. Assume that there
exists a non-trivial automorphism $S$ with a singular spectrum
which is not disjoint from $T$. Then $T$ has a non-trivial factor
which is disjoint from any automorphism with a Lebesgue spectrum.}

The problem of spectral multiplicity of Cartesian products for
``typical" transformation studied  by Katok  \cite{Ka0}  and then
its solution in \cite{Ag3} which we already considered in
Section~\ref{czesc6} lead to a study of those $T$ for which
$$
(CS)\;\;\;\;\sigma^{(m)}\perp\sigma^{(n)}\;\;\mbox{whenever}\;\;m\neq
n,
$$
where $\sigma=\sigma_T$ just stands for the reduced maximal
spectral type of $U_T$ (which is constantly assumed to be a
continuous measure), see also Stepin's article \cite{St}.

Usefulness of the above property (CS) in ergodic theory was
already shown in \cite{Ju-Le}, where a  spectral counterexample
machinery was presented using  the following observation: {\em If
$\ca$ is a $T^{\times\infty}$-invariant sub-$\sigma$-algebra such
that the maximal spectral type on $L^2(\ca)$ is absolutely
continuous with respect to $\sigma_T$ then $\ca$ is contained in
one of the coordinate sub-$\sigma$-algebras $\cb$.} Based on that
in \cite{Ju-Le} it is shown how to construct two weakly isomorphic
actions which are not isomorphic or how to construct two
non-disjoint automorphisms which have no common non-trivial
factors (such constructions were previously known for so called
minimal self-joining automorphisms \cite{Ru5}). See also \cite{Ti}
for extensions of those results to $\Z^d$-actions.

Prikhodko and Ryzhikov  \cite{Pr-Ry} proved that the classical
Chacon transformation enjoys the (CS)~property. The SCS property
defined in the glossary is stronger than  the~(CS) condition
above; the SCS property  implies that the corresponding Gaussian
system $S_{\sigma_T}$ has a simple spectrum. Ageev~\cite{Ag9}
shows that Chacon's transformation satisfies the SCS property;
moreover in \cite{Ag3} he shows that the SCS property is satisfied
generically and he gives a construction of a rank one mixing
SCS-system (see also \cite{Ry200}). In \cite{Le-Pa5} it is proved
that some special flows considered in Section~\ref{czesc12}
(including the von Neumann class, however with $\alpha$ having
unbounded partial quotients) have the SCS property.
It is quite plausible that the SCS property is commonly
seen for smooth flows on surfaces.

A classical open problem is whether each ergodic automorphism has a smooth model. While this problem stays open even for so called dyadic adding machine (ergodic, discrete spectrum automorphism having roots of unity of degree $2^n$, $n\geq1$, as eigenvalues), A.\ Katok suggested many years ago one can
construct a Kronecker measure so that the corresponding Gaussian
system ($\Z$-action (!)) has a smooth representation on the torus. No ``written'' proof of this fact yet appeared.

In \cite{Ve}, A.M.\ Vershik sketches a proof of the fact (claimed by himself for decades) that Pascal adic transformation is weakly mixing. Earlier, X.\ M\'ela in his PhD showed that this transformation is ergodic and has zero entropy. Moreover, in \cite{Ja-Rue}, Janvresse and de la Rue proved the LB property of this map. No complete proof of weak mixing of Pascal adic transformation seems to exist in the published literature.

In \cite{Ve100}, Vershik proposed to study a new equivalence between measure-preserving automorphisms, called  quasi-similarity. Two automorphisms $T$ on $\xbm$ and $S$ on $\ycn$ are {\em quasi-similar} if there are Markov operators $J:L^2\xbm\to L^2\ycn$  and $K:L^2\ycn\to L^2\xbm$, both with dense ranges, intertwining the corresponding Koopman operators $U_T$ and $U_S$. This new equivalence is strictly stronger than spectral isomorphism but strictly weaker than (weak) isomorphism as shown in \cite{Fr-Le100} answering one of questions by Vershik. However, possible invariants for this new equivalence are not well understood. For example, in \cite{Fr-Le100} it is shown that an automorphism quasi-similar to a K-automorphism must also be K, but an intriguing question remains whether Bernoulli property is an invariant of quasi-similarity. As weak isomorphism for finite multiplicity automorphisms is in fact isomorphism, another question is whether we can have two non-isomorphic quasi-similar automorphisms with simple spectrum (rank one)? (see Problem~4 in \cite{Fr-Le100}). See also \cite{Ha-Mo} for the problem of Markov quasi-factors in the class of Abramov automorphisms.

Not too many zero entropy classical dynamical systems with purely Lebesgue spectrum are known: non-zero time automorphisms of horocycle flows and their smooth time changes and even factor of special Gaussian systems. Can we produce a Poissonian suspension over a conservative infinite measure-preserving automorphism having
purely Lebesgue spectrum?

Katok and Thouvenot (private communication) considered systems
called {\em infinitely divisible} (ID). These are systems $T$ on $\xbm$
which have a family of factors  $\cb_{\omega}$ indexed by
$\omega\in \bigcup_{n=0}^\infty\{0,1\}^n$  ($\cb_{\vep}=\cb$) such
that $\cb_{\omega0}\perp\cb_{\omega1},\;\;\cb_{\omega0}\vee
\cb_{\omega1}=\cb_\omega$ and for each $\eta\in\{0,1\}^{\N}$, $
\cap_{n\in\N}\cb_{\eta[0,n]}=\{\emptyset,X\}$.
 They showed (unpublished) that there are discrete spectrum
transformations which are ID, and that there are rank one
transformations with continuous spectra which are also ID (clearly
Gaussian systems are ID). In \cite{Le-Pa-Ro}, it is
proved that dynamical systems coming from stationary ID processes
are factors of ID automorphisms; moreover,  ID automorphisms are
disjoint from all systems having the SCS property. It would be
nice to decide whether Koopman representations associated to ID
automorphisms satisfy the Kolmogorov group property.

While some lacunary sequences have realization as rigidity sequences for weakly mixing automorphisms, it would be interesting to determine what happens in the non-lacunary case. An especially interesting case is when we consider the  non-lacunary multiplicative semigroup $\{2^i3^j:\: i, j \geq 0\}$. Is the corresponding sequence a rigidity sequence? What about $\{2^k+3^\ell:\:k,\ell\geq0\}$ \cite{Be-Ju-Le-Ro}?

\end{document}